%% file: main_v11.tex
\documentclass[a4paper,12pt]{amsart}

\usepackage{amsmath,amsthm,amsfonts,amssymb,mathtools,mathrsfs,bm}
\usepackage{aliascnt}
\usepackage[utf8]{inputenc}
\usepackage[T1]{fontenc}
\usepackage{lmodern}
\usepackage{accents}

\usepackage{tikz-cd}
\usetikzlibrary{babel,decorations.pathmorphing}
\usepackage[shortlabels]{enumitem}
\usepackage{float}
\usepackage{subcaption}
\usepackage{multirow,boldline}
\usepackage{graphicx}
\usepackage{xcolor}

\definecolor{colordelink}{rgb}{0,0,0.50}
\definecolor{colordecite}{rgb}{0,0.5,0}
\definecolor{colordeurl}{rgb}{0,0.41,0.5}
\usepackage[pagebackref=true]{hyperref}
\hypersetup{
  colorlinks=true,
  linkcolor=colordelink,
  citecolor=colordecite,
  urlcolor=colordeurl
}
\usepackage[capitalize,nameinlink,noabbrev]{cleveref}
\usepackage[all]{hypcap}
\usepackage{todonotes} 
\usepackage{adjustbox}

\input{Style-and-more}

\presetkeys{todonotes}{color=yellow!40,linecolor=black!5,size=\small}{}

\DeclareMathOperator{\Crit}{Crit}
\DeclareMathOperator{\Sing}{Sing}
\def\rond{\mathaccent"7017}  

\makeatletter
\@namedef{subjclassname@2020}{\textup{2020} Mathematics Subject Classification}
\makeatother

\title[A general Lê--Siersma formula]{A general Lê--Siersma formula}
\author{R. Giménez Conejero and
J.J.~Nu\~no-Ballesteros}
\address{Escuela Técnica Superior de Ingeniería y Diseño Industrial (ETSIDI),
Universidad Politécnica de Madrid,
Ronda de Valencia, 3
28012 Madrid, Spain
}
\email{roberto.gimenez@uv.es}
\address{Departament de Matemàtiques,
Universitat de Val\`encia, Campus de Burjassot, 46100 Burjassot, Spain}
\email{juan.nuno@uv.es}

\thanks{Work partially supported by Grant PID2025-172404NB-I00 funded by MCIN/AEI/10.13039/501100011033 and by “ERDF A way of making Europe” and by Grant UV-INV\_AE\_4225594}


\begin{document}
\begin{abstract}
We consider the comparison between a Milnor fibre of a function $f$ and
its section by the zero locus of a second function $s$ in a general ambient variety $\mathcal{X}$. We give a formula that generalizes those of Lê and Siersma, expressing the difference of their Euler
characteristics as a sum of local contributions associated with the critical points of the restriction of $s$. We provide an algebraic formula in the case that $\mathcal{X}$ is  {smooth of pure dimension outside $f^{-1}(0)$}, and show some applications and examples.
\end{abstract}

\dedicatory{In memory of Professor L\^e D\~ung Tr\'ang}

\maketitle

\section{Introduction}\label{sec:introduction}

The aim of this paper is to provide a generalization of a result which appears in two different versions in the works of L\^e D\~ung Tr\'ang \cite{Trang1979} and Siersma \cite{Siersma1991}. The first version, due to Lê, is the following:

\begin{theorem} Let $(X,0)$ be the germ of a complete intersection of dimension $n+1$ in $(\C^N,0)$. There exists a Zariski open dense subset $\Omega$ of the dual space $(\C^N)^*$ such that if $\ell\in\Omega$, then the Milnor fibre of the restriction $\ell\colon(X,0)\to(\C,0)$ has the homotopy type of a bouquet of spheres of dimension $n$.
\end{theorem}

In fact, this theorem was generalized later by Hamm and Lê in \cite{Hamm-Le-Handbook} in the case that $X$ is a Milnor space (that is, a space with maximal rectified homotopical depth at any point $x\in X$) and we consider a function $f\colon(X,x)\to(\C,0)$ with an isolated critical point in the stratified sense, instead of a generic linear form (see Theorem \ref{Hamm-Le}).

The second version, due to Siersma, only considers the case that $(X,0)$ is a hypersurface in $(\C^{n+2},0)$, but the function on $(X,0)$ is any linear function with an isolated critical point in the stratified sense. We can assume, without loss of generality, that $\C^{n+2}=\C^{n+1}\times\C$ and the linear function is the projection $\pi\colon \C^{n+1}\times\C\to\C$ given by $\pi(x,t)=t$. Moreover, he gives an explicit formula for the number of spheres of the Milnor fibre of $\pi|_X\colon (X,0)\to(\C,0)$. 

We assume that $s\colon(\C^{n+1}\times\C,0)\to(\C,0)$ is a reduced equation for $(X,0)$. For each $t\in\C$, we put $s_t(x)=s(x,t)$ and $X_t=s_t^{-1}(0)$, so we have a family of hypersurfaces $\{X_t\}$ in $\C^{n+1}$. The Milnor fibre of $\pi|_X$ is $B_\epsilon\cap X_t$, with $0<|t|\ll\epsilon\ll 1$. With this notation, Siersma's result says the following:

\begin{theorem} Let $s\colon(\C^{n+1}\times\C,0)\to(\C,0)$ be a holomorphic function such that $\pi |_X$ has an isolated critical point, where $X=s^{-1}(0)$. Then, $B_\epsilon\cap X_t$ has the homotopy type of a bouquet of spheres of dimension $n$. Moreover, $\beta_n(B_\epsilon\cap X_t)$ is equal to the sum of the Milnor numbers of the critical points of $s_t $ on $B_\epsilon\setminus X_t$. 
\end{theorem}

We remark that the number of spheres of the Milnor fibre appears also implicit in the proof of Lê's result, although he does not mention it in the statement of the theorem. Because of this, we feel that it is fair to call it the Lê--Siersma formula.
Here, we will present a more general formula which includes both results and that we call a general Lê--Siersma formula.
In our case, we consider any germ of an analytic variety $(\Xcal,0)\subseteq(\CC^N,0)$, two functions on it 
$f,s\colon(\Xcal,0)\to(\mathbb C,0)$, and a Whitney stratification of $\Xcal$ adapted to $X:=s^{-1}(0)$. The only assumption is that $f|_X$ has an isolated critical point in the stratified sense. Our main result is
\[
\chi(F_f)-\chi\big(F_{f}\cap X\big)=\sum_i \chi_i n_i
\]
where $n_i$ is the number of stratified Morse critical points in each stratum $S_i$ of a Morsification of the restriction $s_|$ of $s$ to $F_f$ near the finite set $\Gamma\cap F_f$ of critical points of $s_|$, $\Gamma$ being the relative polar curve, and 
$\chi_i$ is, up to the sign, the Euler characteristic of the pair $(\operatorname{Cone}(\mathcal L_i),\mathcal L_i)$, where $\mathcal L_i$ is the complex link of $S_i$  (see Theorem \ref{thm:general}).

In Section 3, we consider the case that $\Xcal$ is a Milnor space and both $f$ and $f|_X$ have an isolated critical point in the stratified sense. In such a case, $F_f$ and $F_f\cap X$ have the homotopy type of a bouquet of spheres. The left hand side of the formula is equal, up to the sign, to the sum of the numbers of spheres of $F_f$ and of $F_f\cap X$. Moreover, in the right hand side, the number $\chi_i$ also counts the number of spheres of $\mathcal L_i$ (see Corollary \ref{cor:Milnor}).

In Section 4, we assume that $\mathcal X\setminus f^{-1}(0)$ is smooth of pure dimension $n+1$. In such case, the general Lê--Siersma formula has an algebraic version as follows:
\[
\chi(F_f)-\chi(F_{f}\cap X)=(-1)^n \dim_\C\frac{\mathcal O_{\Xcal,0}}{(f)+J_{\Xcal}(f,s)\colon (f\cdot s)^\infty} \;,
\]
where $J_{\Xcal}(f,s)$ is the Jacobian ideal of $(f,s)$ in $\Xcal$ (see Theorem \ref{thm:algebraic}). When $f$ has an isolated critical point, the saturation can be taken with respect to $s$ instead of $f\cdot s$ (see Lemma \ref{lem:Gamma}). This formula is very useful, since it can be implemented easily in a computer algebra system like \textsc{Singular}, see \cite{DGPS}. 

Finally, in Section 5, we explain how the algebraic formula can be used to recover the classical Lê--Siersma formula and we also give an application to the image Milnor number $\mu_I(f)$ of a map germ $f\colon(\C^n,S)\to(\C^{n+1},0)$ with isolated instability, as defined by Mond in \cite{Mond1991}.

\section{General result}\label{sec:general}

Let $(\Xcal,0)\subseteq(\CC^N,0)$ be a germ of an analytic variety. We consider {an analytic} Whitney stratification $\Scal=\{S_i\}_i$ of $\Xcal$. Throughout the paper, unless otherwise specified, by a critical point of a holomorphic mapping $h\colon \Xcal\to\C^r$ we will mean a stratified critical point, that is, a critical point of the restriction $h|_{S_i}$, where $S_i$ is the stratum that contains the point. We denote by $\Crit h$ the set of critical points of $h$. Recall that $\Crit h$ is closed, as a consequence of Whitney's condition (a), and analytic.

In this section, we consider two holomorphic functions 
\[
 f,s\colon(\Xcal,0)\longrightarrow(\mathbb C,0).
\]
We will argue with the pair 
\[\Phi=(f,s)\colon(\Xcal,0)\to(\CC^2,0).\]
 To simplify the notation we will also set $X\coloneqq s^{-1}(0)$. We assume that $X$ is a hypersurface in $\Xcal$, i.e., $\dim X=\dim\Xcal -1$, and that the stratification of $\Xcal$ is adapted to $X$, i.e., that $X$ is union of strata. 

Throughout this section we shall also use a Whitney refinement $\mathcal{S}'$ of
$\mathcal{S}$ adapted to $f^{-1}(0)\cup X$ whose strata not contained in
$f^{-1}(0)$ are the sets $S_i\setminus f^{-1}(0)$; it is obtained by
stratifying $S_k\cap f^{-1}(0)$ inside each stratum $S_k$, and the Whitney
conditions for the remaining pairs are inherited from $\mathcal{S}$ (the pairs formed by a stratum of $S_k\cap f^{-1}(0)$ and $S_k\setminus f^{-1}(0)$ satisfy the Whitney conditions because $S_k\setminus f^{-1}(0)$ is open in the manifold $S_k$, and refining the smaller stratum of a pair preserves the Whitney conditions). Since
only the strata $S_i\setminus f^{-1}(0)$ meet $f^{-1}(t)$ for $t\neq 0$, the
induced stratification of $F_f$, the numbers $n_i,\chi_i$ and the polar
curve $\Gamma$ introduced below are the same for $\mathcal{S}$ and $\mathcal{S}'$. Statements
about isolated critical points always refer to $\mathcal{S}$, while the Thom
conditions $a_f$, $a_s$, $a_\Phi$ and the Thom--Mather arguments refer to
$\mathcal{S}'$.

We recall the construction of the Milnor--Lê fibrations of $f$ and $f_|:=f|_X$.
We choose a Milnor radius $0<\varepsilon_0\ll 1$ for $\mathcal S'$, that is, such that every stratum of $\mathcal S'$ is transverse to the sphere $S_\varepsilon$ for all $0<\varepsilon\le\varepsilon_0$; it exists by the curve selection lemma applied to the critical points of $\rho|_{S_j}$, $\rho(x)=\|x\|^2$, for each of the finitely many strata $S_j$ (if
$z(\tau)\to 0$ were an analytic curve of critical points of $\rho|_{S_j}$,
then $\frac{d}{d\tau}\rho(z(\tau))=d\rho(z'(\tau))=0$, so $\rho\circ z$
would be constant, hence zero). In particular, $\varepsilon_0$ is a common Milnor radius for $\Xcal, X$ and $ f^{-1}(0)$. Given $0<\varepsilon\le \varepsilon_0$, we choose a small enough radius $0<\delta_M\ll \varepsilon$ such that, if $t$ is in the punctured disc $D_{\delta_M}^*$, $\Xcal\cap f^{-1}(t)$ is transverse, in the stratified sense, to the Milnor sphere $S_\varepsilon$ and $t$ is a regular value of 
\[f_|\colon(X,0)\to(\CC,0)\]
and, because we will need this technical detail later, of the extension $f^{ext}$ of $f$ to the ambient space $\CC^N$ of $\Xcal$,
\begin{equation}\label{eq:fext}
f^{ext}\colon(\CC^N,0)\to(\CC,0).
\end{equation}
In this situation we have a Milnor--Lê fibration for both $f$ and $f_|$ (see, e.g., \cite[Section 6.11]{Le2020}). Their respective Milnor fibers are
\[
F_f\coloneqq B_\varepsilon\cap\Xcal\cap f^{-1}(t)\quad \text{and} \quad F_{f}\cap X,
\]
for such $t\in D_{\delta_M}^*$. We also put
\[
\rond{F}_f\coloneqq \rond{B}_\varepsilon\cap\Xcal\cap f^{-1}(t)\quad \text{and} \quad \partial F_f\coloneqq S_\varepsilon\cap\Xcal\cap f^{-1}(t).
\]
Since $\Xcal\cap f^{-1}(t)$ is transverse to the spheres $S_r$ for all $r$ close to $\varepsilon$, the open fibre $\rond{F}_f$ is obtained from $F_f$ by removing a collar of the boundary, so both have the same homotopy type; we work with the closed fibre $F_f$ for convenience.

We recall the definition of the Thom condition, which will be necessary to prove our main theorem.

\begin{definition}[see {\cite[II.2.5]{Gibson1976a}}]\label{def:Thom}
Let $\phi\colon N \to N'$ be a smooth map between manifolds $N$ and $N'$ and $X \subset N$ and $X' \subset N'$ two Whitney stratified sets such that $\phi(X) \subset X'$ and the induced map $\phi|\colon X \to X'$ is a stratified map, i.e., submersive onto strata when restricted to strata. We say that $S_i$ is \emph{Thom regular over $S_j$ at $p$ with respect to} $\phi|\colon X \to X'$ if any sequence of points $\{q_n\}_n \subset S_i$ converging to $p \in S_j$ is such that
\begin{equation}\label{eq:thom-condition}
  \ker D_p\, \phi|_{S_j} \subseteq \lim_n \ker D_{q_n}\, \phi|_{S_i}\,,
\end{equation}
when the limit exists. If $\phi|$ is Thom regular for any pair of strata we simply say that the pair of stratifications is a \emph{Thom stratification} of $\phi|$ and that $\phi|$ is a \emph{Thom map} or that it satisfies the \emph{Thom $a_{\phi|}$ condition}.
\end{definition}

It is well known that any function $h\colon(\Xcal,0)\to(\C,0)$ satisfies the Thom condition with respect to any Whitney stratification of $\Xcal$ adapted to $h^{-1}(0)$ (see e.g. \cite[Theorem 4.2.1]{Briancon1994}). However, this is not true in general for mappings $h\colon(\Xcal,0)\to(\C^r,0)$, when $r>1$. There is a classical example due to L\^e (published in a paper of Sabbah \cite{Sabbah1983}) of a mapping $(\C^3,0)\to(\C^2,0)$ which does not admit any stratification which satisfies the Thom condition. The next theorem, due to L\^e, Seade and the second author, shows that if $f|_X$ has an isolated critical point, then it is possible to refine the stratification in such a way that $\Phi=(f,s)\colon (\Xcal,0)\to(\C^2,0)$ satisfies the Thom condition (see \cite[Theorem 5.7]{LNS2025}).

\begin{theorem}\label{thm:LNS} Let  $f,s\colon(\Xcal,0)\longrightarrow(\mathbb C,0)$ be two functions such that $f|_X$ has an isolated critical point. Then:
\begin{enumerate}
\item  $\Gamma\coloneqq\overline{\Crit\Phi\setminus X}$ is empty or a curve in $(\Xcal,0)$;

\item $\Phi^{-1}(0)\cap \Gamma=\{0\}$ and hence,
$\Delta:=\Phi(\Gamma)$
is an analytic curve in $(\C^2,0)$;

\item $\Phi\colon \Xcal\to \C^2$ admits a stratification $\{\mathcal S'',\mathcal T\}$ such that $\mathcal S''$ refines  {$\mathcal S'$} and satisfies the Thom condition.
\end{enumerate}

In particular, $\Phi=(f,s)\colon (\Xcal,0)\to(\C^2,0)$ admits a Milnor--L\^{e} fibration.
\end{theorem}

In Lemma \ref{lem:thom} below we show that refining  {$\mathcal S'$ described above to contain $\Gamma$ as union of strata (and $\Delta$ in the target) }satisfies the Thom condition, which is  {the stratification $\mathcal{S}''$ in which we obtain} Theorem \ref{thm:LNS}(3).

\begin{definition}\label{defn:Gamma} When $f|_X$ has an isolated critical point, the curves $\Gamma$ and $\Delta$ in Theorem \ref{thm:LNS} are called  the \emph{relative polar curve} and  the \emph{Cerf diagram} of $\Phi$, respectively. 
\end{definition}

Note that $\Gamma\cap X=\{0\}$: $\Gamma$ is empty or a curve none of whose irreducible components is contained in $X$, so each component meets $X$ only at the origin. In particular, $\Delta\cap\{v=0\}=\{0\}$, where $(u,v)$ are the coordinates of $\C^2$.

\begin{remark} We remark that the hypothesis that appears originally in \cite[Theorem 5.7]{LNS2025} is that $\Crit(f)\cap X=\{0\}$, instead of $f|_X$ having an isolated critical point. But both conditions are in fact equivalent. This is because the stratification is adapted to $X$ and $\Crit(f)$ is closed (as a consequence of Whitney's condition (a)), hence $\Crit(f)\cap X=\Crit(f|_X)$. 
\end{remark}

We need to introduce a last definition to state and prove our main theorem.

\begin{definition}\label{def:ni-chii}
Suppose that $f|_X$ has an isolated critical point. For each stratum $S_i$ of $\Xcal$, we introduce the following notation:
\begin{enumerate}
\item $n_i$ is the number of stratified Morse critical points in the stratum $S_i$ of a Morsification of the restriction $s_|$ of $s$ to $F_f$ near the finite set $\Gamma\cap F_f$;
\item $\chi_i:=(-1)^{\dim S_i-1}
     \chi\bigl(\operatorname{Cone}(\mathcal L_i),\mathcal L_i\bigr)$, where $\mathcal L_i$ is the complex link of $S_i$ and $\chi(A,B)$ is the Euler characteristic of the pair $(A,B)$. When $S_i$ is open in $\Xcal$ we have $\mathcal L_i=\varnothing$ and $\operatorname{Cone}(\varnothing)$ is a point, so that $\chi_i=(-1)^{\dim S_i-1}$ in this case.
\end{enumerate}
\end{definition}

\begin{theorem}\label{thm:general}
Let  $f,s\colon(\Xcal,0)\longrightarrow(\mathbb C,0)$ be two functions such that $f|_X$ has an isolated critical point. Then,
\begin{equation}\label{eq:main-formula}
 \chi(F_f)-\chi\big(F_{f}\cap X\big)=\sum_i \chi_i n_i.
\end{equation}
\end{theorem}

The proof of Theorem \ref{thm:general} is based on stratified Morse theory. We apply it to a map  {induced by the restriction of $|s|$ to $F_f$}.
We need several technical lemmas to ensure that the critical points on the boundary do not contribute to the change of the topology of the fibre.

\begin{lemma}\label{lem:crit-s}
For every stratum $S_j$ of $\mathcal S'$ contained in $f^{-1}(0)$, the germ at $0$ of $\Crit(s|_{S_j})$ is contained in $X$.
\end{lemma}

\begin{proof}
Let $C$ be an irreducible component of the closed analytic set germ $\Crit(s)\cap f^{-1}(0)$ (the critical locus being taken with respect to $\mathcal S'$). Then $C$ passes through $0$, and there is a stratum $S_k$ such that $C\cap S_k$ is open and dense in $C$. Since $d(s|_{S_k})=0$ on $C\cap S_k$, which is connected (it is the complement of a proper analytic subset of the irreducible germ $C$), $s$ is constant on $C$, and therefore $s\equiv s(0)=0$ on $C$. Thus $\Crit(s)\cap f^{-1}(0)\subseteq X$, and in particular $\Crit(s|_{S_j})\subseteq X$.
\end{proof}

\begin{lemma}\label{lem:thom}
Assume that $f|_X$ has an isolated critical point in the stratified sense (with respect to $\mathcal S$).
Then $\Phi=(f,s)$ satisfies the Thom $a_\Phi$ condition with respect to
$\mathcal{S}'$ at every point of $f^{-1}(0)\cup X$. 
\end{lemma}

\begin{proof}
Let $S'_j\subseteq f^{-1}(0)\cup X$ and $S'_i$ be strata of $\mathcal S'$ with $S'_j\subseteq\overline{S'_i}$, let $p\in S'_j$ and let $q_n\in S'_i$ be a sequence
converging to $p$. We may assume that
$p\neq 0$, since $\ker D_0\Phi|_{\{0\}}=0$. We denote by $S_j\supseteq S'_j$ and $S_i\supseteq S'_i$ the strata of $\mathcal S$ containing them. Passing to a subsequence, we
assume that $T_{q_n}S_i\to T$, so that $T\supseteq T_pS_j$ by the Whitney
condition (a), and that
\[
  A_n:=\ker D_{q_n}(s|_{S_i})\to A,\qquad
  H_n:=\ker D_{q_n}(f|_{S_i})\to H,
\]
where we also assume, passing to a subsequence, that either $A_n$ is a
hyperplane of $T_{q_n}S_i$ for all $n$ or $A_n=T_{q_n}S_i$ for all $n$, and
similarly for $H_n$. 
It is not difficult to check that
\begin{equation}\label{eq:intersection}
A\ne H\Longrightarrow A_n\cap H_n\to A\cap H.
\end{equation}
We distinguish two cases, according to whether $S'_i$ is contained in $f^{-1}(0)$ or not.

\emph{Case $S'_i\not\subseteq f^{-1}(0)$.} Then $S'_i=S_i\setminus f^{-1}(0)$ is open in $S_i$, so that $\ker D_{q_n}\Phi|_{S'_i}=\ker D_{q_n}\Phi|_{S_i}=A_n\cap H_n$. Note that if $A_n=T_{q_n}S_i$ for all $n$ then $\lim_n\ker D_{q_n}\Phi|_{S_i}=H$, and symmetrically. We consider two subcases.

\emph{Subcase $p\in X$.} Then $S_j\subseteq X$ and $\ker D_p\Phi|_{S'_j}\subseteq T_pS_j\cap\ker d_pf$, with equality if $S'_j=S_j\setminus f^{-1}(0)$ (if $S'_j\subseteq f^{-1}(0)$, then $\ker D_p\Phi|_{S'_j}\subseteq T_pS'_j\subseteq T_pS_j\cap\ker d_pf$). Since $\mathcal{S}$ is adapted to
$X=s^{-1}(0)$, the Thom condition $a_s$ holds for $\mathcal S$ by \cite[Theorem
4.2.1]{Briancon1994}, so $A\supseteq T_pS_j$. Since $f|_{S_j}$ is a submersion at
$p$ (because $p\ne0$ and $f|_X$ has an isolated critical point) and $T\supseteq T_pS_j$, we have $d_pf|_T\neq 0$, hence $H_n$ are
hyperplanes for $n\gg 0$, $H=\ker(d_pf|_T)$ and $H\not\supseteq T_pS_j$.
If $A_n=T_{q_n}S_i$ for all $n$, then
$\lim_n\ker D_{q_n}\Phi|_{S_i}=H\supseteq T_pS_j\cap\ker d_pf
\supseteq\ker D_p\Phi|_{S'_j}$. Otherwise $A\neq H$ and 
\eqref{eq:intersection} gives
\[
  \lim_n\ker D_{q_n}\Phi|_{S_i}=A\cap H\supseteq T_pS_j\cap\ker d_pf
  \supseteq\ker D_p\Phi|_{S'_j}.
\]

\emph{Subcase $p\notin X$.} Then $S'_j\subseteq f^{-1}(0)\setminus X$ and $\ker D_p\Phi|_{S'_j}=T_pS'_j\cap \ker d_ps$, because $f|_{S'_j}\equiv0$. Since $\mathcal{S}'$ is
adapted to $f^{-1}(0)$, the Thom condition $a_f$ holds for $\mathcal S'$ by \cite[Theorem
4.2.1]{Briancon1994}, so $H\supseteq T_pS'_j$ (in particular $H_n$ are hyperplanes). 
By Lemma \ref{lem:crit-s}, $d_p(s|_{S'_j})\neq 0$; since $T\supseteq T_pS_j\supseteq T_pS'_j$, we get $d_ps|_T\neq 0$, so that $A_n$ are hyperplanes for $n\gg0$, $A=\ker(d_ps|_T)$ and
$A\not\supseteq T_pS'_j$. Hence $A\neq H$ and \eqref{eq:intersection} gives
\[
  \lim_n\ker D_{q_n}\Phi|_{S_i}= A\cap H\supseteq T_pS'_j\cap\ker d_ps
  =\ker D_p\Phi|_{S'_j}.
\]

\emph{Case $S'_i\subseteq f^{-1}(0)$.} Then also $S'_j\subseteq f^{-1}(0)$, so $f$ vanishes identically on both strata and $\ker D_{q_n}\Phi|_{S'_i}=\ker D_{q_n}(s|_{S'_i})$, $\ker D_p\Phi|_{S'_j}=\ker D_p(s|_{S'_j})$. The required condition is therefore the Thom condition $a_s$ for the restriction of $s$ to the analytic space $f^{-1}(0)\cap\Xcal$, endowed with the Whitney stratification induced by $\mathcal S'$, which is adapted to $(f^{-1}(0)\cap\Xcal)\cap s^{-1}(0)=\Phi^{-1}(0)$. If $S'_j\subseteq X$, it holds by \cite[Theorem 4.2.1]{Briancon1994} applied to this space. If $S'_j\not\subseteq X$, then $d_p(s|_{S'_j})\ne0$ by Lemma \ref{lem:crit-s}, hence $d_ps$ does not vanish on $T':=\lim_n T_{q_n}S'_i\supseteq T_pS'_j$ (Whitney condition (a) for $\mathcal S'$), and $\lim_n \ker D_{q_n}(s|_{S'_i})=\ker(d_ps|_{T'})\supseteq \ker D_p(s|_{S'_j})$.
\end{proof}

For a stratum $S_i$ of $\mathcal S'$ and
$q\in S_i$ at which $\Phi|_{S_i}$ is a submersion, the fibre
$\Phi^{-1}(\Phi(q))\cap S_i$ is a smooth manifold near $q$ with tangent space
$\ker D_q\Phi|_{S_i}$.  {In this case, we} say that \textit{the fibres of $\Phi$ are transverse to
$S_\varepsilon$ over a subset} $W\subseteq\mathbb{C}^2$ if
\begin{equation}\label{eq:transv}
  \ker D_q\Phi|_{S_i}\not\subseteq T_qS_\varepsilon
\end{equation}
for every stratum $S_i$ and every $q\in S_i\cap\Phi^{-1}(W)\cap S_\varepsilon$
at which $\Phi|_{S_i}$ is a submersion.

\begin{lemma}\label{lem:C1}
Let $\varepsilon_0>0$ be the Milnor radius for $\mathcal S'$ fixed above. For every $0<\varepsilon\le\varepsilon_0$ there is
$\delta(\varepsilon)>0$ such that the fibres of $\Phi$
are transverse to $S_\varepsilon$ over the polydisc
$D_{\delta(\varepsilon)}\times D_{\delta(\varepsilon)}$  {(as defined above, at strata $S_i$ where $\Phi|_{S_i}$ is a submersion)}.
\end{lemma}

\begin{proof}
Fix $0<\varepsilon\le\varepsilon_0$ and suppose that no $\delta(\varepsilon)$
exists. Then there are a stratum $S_i$ of $\mathcal S'$ (fixed after passing to a
subsequence), points $q_n\in S_i\cap S_\varepsilon$ at which $\Phi|_{S_i}$ is
a submersion, with $\Phi(q_n)\to 0$ and
$\ker D_{q_n}\Phi|_{S_i}\subseteq T_{q_n}S_\varepsilon$. By compactness we
may assume that $q_n\to p_0\in S_\varepsilon$ and that
$\ker D_{q_n}\Phi|_{S_i}\to K$ in the corresponding Grassmannian. Then
$\Phi(p_0)=0$, so $p_0$ lies in a stratum $S_j\subseteq\Phi^{-1}(0)$ of $\mathcal S'$ with
$S_j\subseteq\overline{S_i}$, and $K\subseteq T_{p_0}S_\varepsilon$ since
$T_{q_n}S_\varepsilon\to T_{p_0}S_\varepsilon$. On the other hand, the Thom
condition $a_\Phi$ for the pair $(S_j,S_i)$ at $p_0$ (Lemma \ref{lem:thom}) gives
\[
  K\supseteq\ker D_{p_0}\Phi|_{S_j}=T_{p_0}S_j ,
\]
because $\Phi|_{S_j}\equiv 0$. Hence $T_{p_0}S_j\subseteq T_{p_0}S_\varepsilon$,
contradicting the transversality of $S_j$ to $S_\varepsilon$.
\end{proof}

\begin{remark}\label{rem:C1-consequences}
Let $\varepsilon\le\varepsilon_0$ and $|t|,|v|\le\delta(\varepsilon)$. The
following consequences of Lemma \ref{lem:C1} are the ones used below.
\begin{enumerate}
\item If $f|_{S_i}$ is a submersion at $q\in S_i\cap f^{-1}(t)\cap S_\varepsilon$
  and $s(q)=v$ with $d_q(s|_{S_i\cap f^{-1}(t)})\neq 0$, then the level set of
  $|s|$ in $S_i\cap f^{-1}(t)$ through $q$ is transverse to $S_\varepsilon$ at
  $q$, since its tangent space contains $\ker D_q\Phi|_{S_i}$. In particular,
  $|s|$ has no critical points on the boundary strata
  $S_i\cap f^{-1}(t)\cap S_\varepsilon$ of $F_f$ in the region
  $\{|s|\le\delta(\varepsilon)\}$, except possibly at points of
  $\Gamma\cap f^{-1}(t)$, which do not meet $S_\varepsilon$ for $|t|$ small
  because $\Phi|_\Gamma$ is finite.
\item The restriction of $s$ to each boundary stratum
  $S_i\cap f^{-1}(t)\cap S_\varepsilon$ of $F_f$ is a submersion onto
  $D_{\delta(\varepsilon)}$ away from $X$, and the Milnor fibre
  $F_f\cap X=\Phi^{-1}(t,0)\cap B_\varepsilon$ of $f|_X$ is
  transverse to $S_\varepsilon$.
\end{enumerate}
\end{remark}

We regard $F_f=\mathcal{X}\cap f^{-1}(t)\cap B_\varepsilon$ as a
Whitney stratified set with strata
\[
  S_i^{\circ}:=S_i\cap f^{-1}(t)\cap \rond{B}_\varepsilon
  \qquad\text{and}\qquad
  S_i^{\partial}:=S_i\cap f^{-1}(t)\cap S_\varepsilon ,
\]
where $S_i$ runs over the strata of $\mathcal S'$ not contained in $f^{-1}(0)$ (equivalently, over the strata of $\mathcal S$);
this is a Whitney stratification because $t$ is a regular value of $f$
and $f^{\mathrm{ext}}$, and $\Xcal\cap f^{-1}(t)$ is transverse to the Milnor sphere
(see \eqref{eq:fext}). We denote by $g$ the restriction of $|s|$ to $F_f$.

We take a point $x\in F_f\setminus X$. We want to characterize the critical points of $g$. 
When $x\in S_i^\circ$, it is well known that $x$ is a critical point of $g|_{S_i^\circ}$ if and only if $x$ is a critical point of $s|_{S_i^\circ}$. Otherwise, if $x\in S_i^\partial$, we say that $x$ is a \emph{boundary critical point} of $g$ if it is a critical point of the restriction $g|_{S_i^\partial}$. We also set 
\[
K_x\coloneqq T_x\big(S_i\cap f^{-1}(t)\big),
\]
so that $T_xS_i^\partial=K_x\cap T_xS_\varepsilon$.

\begin{lemma}\label{lem:lagrange} With the above notation, $x\in S_i^\partial$ is a boundary critical point of $g$ if and only if
\begin{equation}\label{eq:lagrange}
  \bar s\,ds|_{K_x}=\lambda\,\partial\rho|_{K_x},
  \end{equation}
for some $\lambda\in\R$  {and $\rho=\|\bullet\|^2$}.
\end{lemma}

\begin{proof} Since $S_i^\partial$ is the level set $\{\rho=\varepsilon^2\}$ of $\rho$ in $S_i\cap f^{-1}(t)$, the point $x$ is a critical point of $g|_{S_i^\partial}$ if and only if
\[
d|s|\,|_{K_x}=\mu\, d\rho|_{K_x},
\]
for some $\mu\in\R$. But $d|s|=\operatorname{Re}(\bar s\,ds)/|s|$ and $d\rho=2\operatorname{Re}(\partial\rho)$, and a $\C$-linear form on a complex vector space is determined by its real part. Hence the above equality is equivalent to
\[
  \bar s\,ds|_{K_x}=\lambda\,\partial\rho|_{K_x},
  \]
where  $\lambda=2\mu |s(x)|$.
\end{proof}

A boundary critical point $x\in\partial{F}_f$ of $g$ is called a \emph{tangency point}. It is called an \emph{outward tangency point} if \eqref{eq:lagrange} holds with $\lambda>0$ (equivalently $\mu>0$), that is, if $|s|$ increases along the outward normal direction of $S_\varepsilon$.

\begin{lemma}\label{lem:outward}
There is $\varepsilon_1>0$ such that for all $\varepsilon\le\varepsilon_1$
and all $\eta_0>0$ there is $\delta=\delta(\varepsilon,\eta_0)>0$ such that if $0<|t|<\delta$, then all tangency points $x\in \partial F_f$ with 
$|s(x)|\ge\eta_0$ are outward tangency points.
\end{lemma}

\begin{proof}
First, $\lambda\neq 0$ at any tangency point for $|t|$ small: otherwise $ds|_{K_x}=0$ and
$x\in\Crit\Phi \setminus X\subseteq\Gamma$, 
which does not meet $S_\varepsilon$ for $|t|$ small because $\Phi|_\Gamma$ is
finite.

We claim that we can choose $\varepsilon_1\le\varepsilon_0$ such that, for every stratum
$S_j$ of $\mathcal S'$ with $S_j\subseteq f^{-1}(0)\setminus X$: 
\begin{enumerate}
\item $\Crit(s|_{S_j})\cap
B_{\varepsilon_1}=\varnothing$; and
\item the semianalytic set
\[
B_j=\{q\in S_j:\ \bar s\,ds|_{T_qS_j}=\lambda\,\partial\rho|_{T_qS_j},\
\lambda<0\}
\] 
does not meet $B_{\varepsilon_1}$. 
\end{enumerate}
Item (1) is Lemma \ref{lem:crit-s}, since $S_j\cap X=\varnothing$. Alternatively, if (1) was false, by the curve selection lemma we could find an analytic curve $\gamma(\tau)$ such that $\gamma(0)=0$, $\gamma(\tau)\in S_j$ and $\gamma(\tau)\in \Crit(s|_{S_j})$, for all $\tau\ne0$; then $\frac{d}{d\tau}s(\gamma(\tau))=d_{\gamma(\tau)}s(\gamma'(\tau))=0$, so $s\circ\gamma$ would be constant, hence zero, and $\gamma(\tau)\in X$, in contradiction with $S_j\cap X=\varnothing$.

For (2), we use again the curve selection lemma. If (2) was false, we could find an analytic curve
$\gamma(\tau)$ such that $\gamma(0)=0$ and $\gamma(\tau)\in B_j$, for all $\tau\ne0$. Taking real parts in the relation defining $B_j$ gives
$\frac{d}{d\tau}|s(\gamma(\tau))|^2=\lambda(\tau)\frac{d}{d\tau}\rho(\gamma(\tau))$; both $|s\circ\gamma|^2$ and $\rho\circ\gamma$ are real analytic, non-negative, not identically zero and vanish at $\tau=0$, so both
derivatives are positive for small $\tau>0$, and hence $\lambda(\tau)>0$, a contradiction.

Suppose now, by contradiction, that the lemma fails for some $\varepsilon\le\varepsilon_1$ and some $\eta_0$. 
There would be a stratum $S_i$ (fixed after passing to a subsequence) and points $x_n\in S_i\cap f^{-1}(t_n)\cap S_\varepsilon$, $t_n\to 0$, with $\lambda_n<0$ and $|s(x_n)|\ge\eta_0$. 

Passing to a subsequence, $x_n\to x\in S_\varepsilon$ and $K_{x_n}\to K$. Then $f(x)=0$ and $|s(x)|\ge\eta_0>0$, so $x$ lies in a stratum $S_j$ of $\mathcal S'$ with $S_j\subseteq f^{-1}(0)\setminus X$. By the
$a_f$ condition for $\mathcal S'$ \cite[Theorem 4.2.1]{Briancon1994}, $K\supseteq T_{x}S_j$. Since
$\partial\rho|_{T_{x}S_j}\neq 0$ (transversality of $S_j$ to
$S_\varepsilon$), we can choose $v_n\in K_{x_n}$ converging to some $v\in T_xS_j$ with $\partial\rho(v)\ne0$, and evaluating \eqref{eq:lagrange} at $v_n$ shows that the $\lambda_n$ are bounded. Passing to a subsequence, $\lambda_n\to\lambda\le0$ and in the limit
$\bar s\,ds=\lambda\,\partial\rho$ on $T_{x}S_j$. If $\lambda<0$, then $x\in B_j\cap B_{\varepsilon_1}$,
contradicting (2); if $\lambda=0$, then $d_x(s|_{S_j})=0$ (recall $s(x)\ne0$) and
$x\in\Crit(s|_{S_j})\cap B_{\varepsilon_1}$, contradicting (1).
\end{proof}

\begin{figure}[htbp]
\centering
\begin{tikzpicture}[scale=.85,line cap=round,line join=round,font=\small,
  Vfield/.style={-stealth,thick,draw=black!75},
  grad/.style={-stealth,thick,draw={rgb,255:red,30;green,140;blue,30}},
  bad/.style={-stealth,very thick,draw={rgb,255:red,190;green,40;blue,40}},
  level/.style={thick,draw=blue!60!black}]
\colorlet{ncDark}{blue!40}
\colorlet{ncLight}{blue!13}
\begin{scope}
  \coordinate (c) at (-0.7,0); 
  \coordinate (p) at (-2.4,0); 
  \begin{scope}
    \clip (0,0) circle (2.4);
    \fill[ncLight] (c) circle (2.45);
    \draw[level] (c) circle (2.45);                       
    \draw[densely dashed,blue!60!black] (c) circle (1.7); 
  \end{scope}
  \fill[ncDark] (c) circle (1.05);
  \draw[level] (c) circle (1.05);                         
  \draw[very thick] (0,0) circle (2.4);
  \foreach \ang/\r in {40/2.05,-40/2.05,80/2,-80/2}
    \draw[Vfield] (c) ++(\ang:\r) -- ++(\ang+180:.5);
  \draw[Vfield] (112:2.4) arc (112:130:2.4);
  \draw[Vfield] (147:2.4) arc (147:166:2.4);
  \draw[Vfield] (248:2.4) arc (248:230:2.4);
  \draw[Vfield] (213:2.4) arc (213:194:2.4);
  \draw[Vfield] (p) -- ++(.6,0);
  \draw[grad] (p) -- ++(-.8,0) node[pos=1,above] {$\nabla g$};
  \fill (p) circle (1.8pt) node[below left] {$p$};
  \fill (c) circle (1.6pt) node[below] {$F_f\cap X$};
  \node at (-0.7,0.45) {$g^{-1}[0,a]$};
  \draw[thin] (1.1,0.6) -- (2.55,1.55) node[right] {$g^{-1}[a,b]$};
  \node[right] at (1.75,-1.95) {$\partial F_f$};
  \node at (0,-3) {(a) outward tangency, $\mu>0$};
\end{scope}
\begin{scope}[xshift=7.4cm]
  \coordinate (c) at (0.7,0);   
  \coordinate (o) at (-3.6,0);  
  \coordinate (p) at (-2.4,0);  
  \begin{scope}
    \clip (0,0) circle (2.4);
    \fill[ncLight] (c) circle (1.35);
    \fill[ncLight] (o) circle (1.7);
    \draw[level] (o) circle (1.7);
  \end{scope}
  \fill[ncDark] (c) circle (.95);
  \draw[level] (c) circle (1.35);
  \draw[level] (c) circle (.95);
  \draw[very thick] (0,0) circle (2.4);
  \draw[grad] (p) -- ++(.75,0) node[right] {$\nabla g$};
  \fill (p) circle (1.8pt) node[below left] {$p$};
  \fill (c) circle (1.6pt) node[below] {$F_f\cap X$};
  \node at (0.7,0.45) {$g^{-1}[0,a]$};
  \draw[thin] (-2.28,0.55) -- (-3.0,2.2)
     node[above] {new piece of $g^{-1}[0,b]$};
  \node[right] at (1.75,-1.95) {$\partial F_f$};
  \node at (0,-3) {(b) inward tangency, $\mu<0$};
\end{scope}
\end{tikzpicture}
\caption{The situation of Lemma~\ref{lem:no-contribution}, with $g=g$ on
the Milnor fibre $F_f$, drawn schematically as a disc. (a) The level set of
$g$ through an outward tangency point $p\in\partial F_f$ (dashed) touches
$\partial F_f$ from inside, and $g$ increases in the outward normal
direction. The vector field $V$ of the proof decreases $g$ at unit speed:
it is tangent to $\partial F_f$ at the other boundary points and points into
$F_f$ at $p$, so its flow retracts $g^{-1}[0,b]$ onto $g^{-1}[0,a]$.
(b) At an inward tangency point, $g^{-1}[0,b]$ acquires a new piece at the
boundary, and the only direction decreasing $g$ at $p$ leaves $F_f$. This
case is excluded by hypothesis~(2). Black arrows: $V$; green arrows: the
gradient of $g$ along $F_f$.}
\label{fig:no-contribution}
\end{figure}
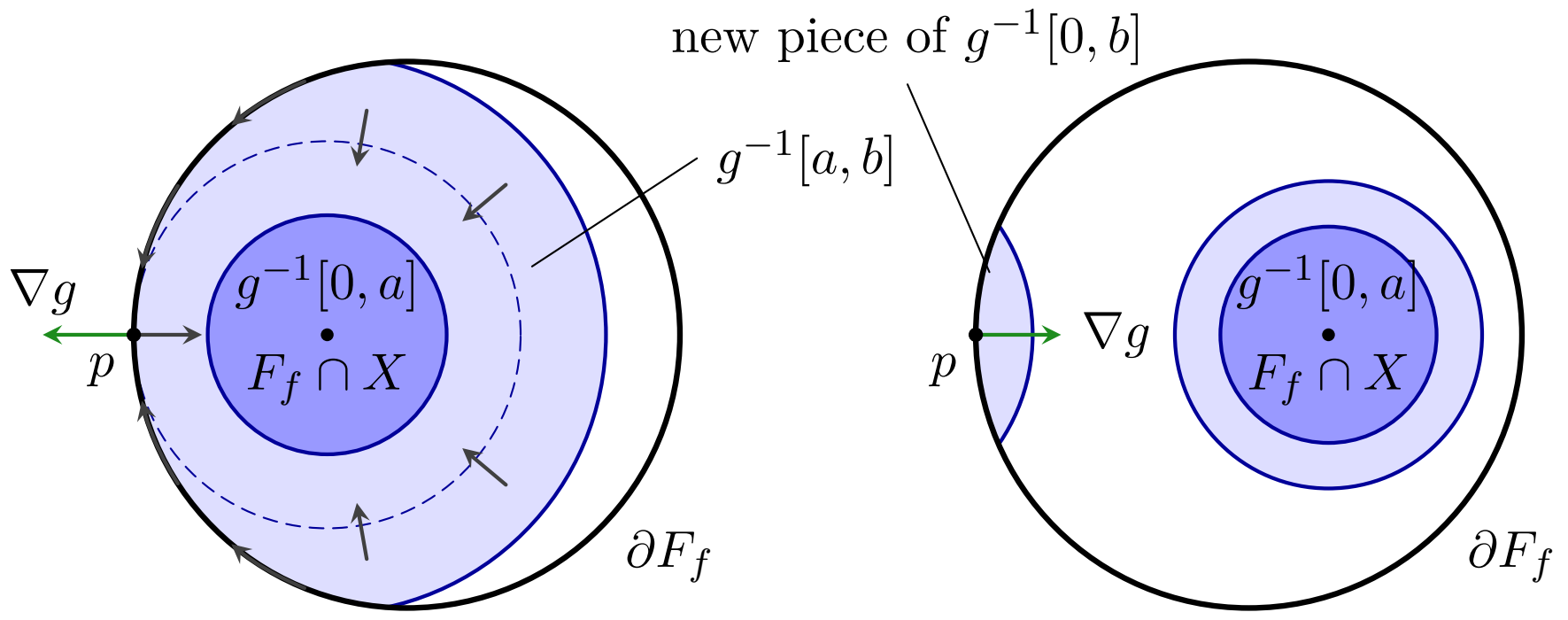

\begin{lemma}\label{lem:no-contribution}
Let $0<a<b$ be such that
\begin{enumerate}
\item $g$ has no critical points on $S_i^{\circ}\cap g^{-1}[a,b]$, for every
  stratum $S_i$;
\item every tangency point $p\in S_i^{\partial}\cap g^{-1}[a,b]$
  is an outward tangency point.
\end{enumerate}
Then $g^{-1}[0,a]$ is a deformation retract of $g^{-1}[0,b]$ (see \cref{fig:no-contribution}). In particular,
the inclusion $g^{-1}[0,a]\hookrightarrow g^{-1}[0,b]$ is a homotopy
equivalence.
\end{lemma}

\begin{proof}
It suffices to construct, on a neighbourhood of $g^{-1}[a,b]$ in $F_f$, a
controlled stratified vector field $V$ in the sense of
\cite[II.4]{Gibson1976a}, tangent to the strata $S_i^{\circ}$ and
$S_i^{\partial}$, such that
\begin{equation}\label{eq:V}
  V( {g})=-1 \qquad\text{and}\qquad V(\rho)\le 0 .
\end{equation}
Indeed, the flow $\varphi_\tau$ of $V$ is continuous, preserves the strata
and, by \eqref{eq:V}, decreases $ {g}$ at unit speed without increasing
$\rho$; hence it does not leave $B_\varepsilon$ and is defined on
$g^{-1}[a,b]$ for $0\le\tau\le b-a$. Then
$(x,\tau)\mapsto\varphi_{\min\{\tau,\,g(x)-a\}}(x)$ is the required
deformation retraction. Since both conditions in \eqref{eq:V} are preserved
by convex combinations, and controlled vector fields with respect to a
fixed system of control data form a module over the smooth functions, it is
enough to construct $V$ locally and glue with a partition of unity. Let
$p\in g^{-1}[a,b]$.

\emph{Case 1: $p\in S_i^{\circ}$.} By (1) and the Whitney condition (a),
$g$ is a stratified submersion in a neighbourhood of $p$, and the Thom--Mather
lifting lemma gives a controlled lift $V$ of the vector field $-d/dr$ on $\R$, $r$ being the coordinate of $\R$; it
satisfies $V( {g})=-1$, and the condition on $\rho$ is  {unnecessary} since $p$ is interior
(shrink the neighborhood so that it does not meet $S_\varepsilon$).

\emph{Case 2: $p\in S_i^{\partial}$ is not a tangency point.} Then
$d_p {g}$ and $d_p\rho$ are linearly independent on $K_p$, so that
$( {g},\rho)\colon F_f\to\mathbb{R}^2$ is a stratified submersion near $p$:
on $S_i^{\partial}=\{\rho=\varepsilon^2\}\cap (S_i\cap f^{-1}(t))$ because
$d_p {g}\,|_{T_pS_i^{\partial}}\neq 0$, on $S_i^{\circ}$ by continuity, and on
the larger strata by the Whitney condition (a) and openness. A controlled
lift $V$ of the constant vector field $(-1,0)$ on $\mathbb{R}^2$ satisfies
$V( {g})=-1$ and $V(\rho)=0$.

\emph{Case 3: $p\in S_i^{\partial}$ is an outward tangency point.} Since all
the strata of $f^{-1}(t)\cap\mathcal{X}$ are transverse to the spheres $S_r$
for $r$ close to $\varepsilon$, the function $\rho$ is a stratified
submersion near $p$, and there is a controlled lift $\xi$ of
$\partial/\partial\rho$, tangent to the strata and continuous. At $p$ we have
$\xi_p\in K_p$ with $d_p\rho(\xi_p)=1$, hence
\[
  \xi( {g})(p)=d_p {g}(\xi_p)=\mu\,d_p\rho(\xi_p)=\mu>0 .
\]
As $\xi$ is continuous and $ {g}$ is smooth near $p$, $\xi( {g})>0$ on a
neighbourhood $U_p$ of $p$, on all strata. Then $V:=-\xi/\xi( {g})$ on $U_p$
satisfies $V( {g})=-1$ and $V(\rho)=-1/\xi( {g})<0$.
\end{proof}

\begin{remark}\label{rem:inward}
Observe that only the sign of $\mu$ is used in Lemma \ref{lem:no-contribution}, it does not require the tangency points on the
boundary to be isolated or non-degenerate.

With $\mu<0$ the vector field decreasing $|s|$ would point outwards, see \cref{fig:no-contribution} (b), and
the sublevel sets would acquire cells at the boundary. This is what happens when the hypothesis of Theorem \ref{thm:general} fails. For instance, take $\Xcal=\C^2$, $f=x^2$ and $s=x(1+z^2)$, so that $X=\{x=0\}$ and $f|_X\equiv0$. Then $F_f=\{x=\pm\sqrt t\}\cap B_\varepsilon$ is the disjoint union of two discs and $F_f\cap X=\varnothing$, so the left hand side of \eqref{eq:main-formula} equals $2$. On the other hand, $\Gamma=\{z=0\}$ meets $F_f$ in two points, which are Morse points of $s|_{F_f}$, and $\chi_0=(-1)^{2-1}=-1$ for the open stratum, so the right hand side equals $-2$. The discrepancy comes from the boundary: on each disc, $|s|=|t|^{1/2}|1+z^2|$ has two minima on the boundary circle which are inward tangency points.
\end{remark}

\begin{corollary}\label{cor:polydisc-fibre}
With $\varepsilon\le\varepsilon_1$ as in Lemma \ref{lem:outward} and
$\delta'\le\delta(\varepsilon)$ as in Lemma \ref{lem:C1}, there exists $\delta_\Gamma=\delta_\Gamma(\varepsilon,\delta')>0$ such that for all 
$0<|t|<\min\{\delta_M,\delta(\varepsilon,\delta'),\delta_\Gamma\}$, the inclusion
\[
  F'_f:=F_f\cap\{|s|\le\delta'\}\hookrightarrow F_f
\]
is a homotopy equivalence.
\end{corollary}

\begin{proof} Since $\Phi|_\Gamma$ is finite, there exists $\delta_\Gamma=\delta_\Gamma(\varepsilon,\delta')>0$ such that for all $0<|t|<\delta_\Gamma$,
\[
\Gamma\cap f^{-1}(t) \cap B_\varepsilon\subset \{|s|<\delta'\}.
\]
Apply Lemma \ref{lem:no-contribution} with $a=\delta'$ and
$b=\max_{F_f}|s|$. Condition (1) holds because the critical points of $g$ on
the interior strata are the points of $\Gamma\cap F_f$, which satisfy
$|s|<\delta'$, for $|t|<\delta_\Gamma$. Condition (2)
is Lemma \ref{lem:outward} with $\eta_0=\delta'$.
\end{proof}

\begin{proof}[Proof of Theorem \ref{thm:general}]
We denote the coordinates of $(\CC^2,0)$ by $(u,v)$. We first fix the constants. We take $\varepsilon\le \min\{\varepsilon_0,\varepsilon_1\}$, where $\varepsilon_0$ is the Milnor radius for $\mathcal S'$ and $\varepsilon_1$ is given by Lemma \ref{lem:outward}, and we let $\delta_M$ be the radius of the Milnor--Lê fibrations of $f$ and $f_|$ in $B_\varepsilon$, as at the beginning of the section. Next we fix 
$\delta'\le\delta(\varepsilon)$ as in Lemma \ref{lem:C1}, and finally we fix $\delta$ such that 
\[
0<\delta<\min\{\delta_M,\ \delta(\varepsilon),\ \delta(\varepsilon,\delta'),\ \delta_\Gamma(\varepsilon,\delta')\},
\]
where $\delta(\varepsilon,\delta')$ and $\delta_\Gamma(\varepsilon,\delta')$ are given by Lemma \ref{lem:outward} and Corollary \ref{cor:polydisc-fibre}. Note that the order of the choices is essential, since the last two constants depend on $\delta'$. We choose and fix a $t_0$ in $D_\delta^*$ and we set
\[
F'_f:=\Phi^{-1}\big(\{t_0\}\times D_{\delta'}\big)\cap B_\varepsilon=F_f\cap\{|s|\le\delta'\}.
\]
By Corollary \ref{cor:polydisc-fibre}, the inclusion $F'_f\hookrightarrow F_f$ is a homotopy equivalence, and obviously $F'_f\cap X=F_f\cap X$ (see \cref{fig:sketch} below).

Regarding the fibre $F_f\cap X$, we can choose a small enough radius $0<\eta<\delta'$ so that the disc $D_\eta$ centered at the origin of the line $\{u=t_0\}$ avoids $\Delta$. This is possible since $(t_0,0)\notin\Delta$: as observed after Definition \ref{defn:Gamma}, $\Gamma\cap X=\{0\}$, so $\Delta\cap\{v=0\}=\{0\}$. In addition, the fibres of $\Phi$ over $\{t_0\}\times D_{\delta'}$ are transverse to the Milnor sphere $S_\varepsilon$ by Lemma \ref{lem:C1}. We show now that 
\[\Phi^{-1}\big(\{t_0\}\times D_\eta\big)\cap B_\varepsilon\]
 is homotopy equivalent to the fiber $F_f\cap X$.

First, observe that the restriction 
\[\Phi_|:\Phi^{-1}\big(\{t_0\}\times D_\eta\big)\cap B_\varepsilon\to \{t_0\}\times D_\eta\]
coincides with the restriction
\[s_|:s^{-1}(D_\eta)\cap F_f\to D_\eta.\]
We will consider the radial vector field on $D_\eta$, whose flow retracts the disc to the origin, and we want to lift it by $s_|$ to a globally integrable vector field on $s_|^{-1}(D_\eta)\cap F_f$ whose flow will retract it to the zero fiber
\[s_|^{-1}(0)=F_f\cap X.\]

 We can lift the vector field if the map $s_|$ is a Thom proper map with certain stratifications on source and target. One can see this in \cite[II.2.6, II.3.2, II.4.6 and II.4.8]{Gibson1976a} or in the more detailed account \cite[Corollary 3.11]{GimenezConejero2023}.

To see that the map $s_|$ is proper, observe that if $K\subset D_\eta$ is compact, then $s_|^{-1}(K)=F_f\cap s^{-1}(K)$ is closed in the closed bounded set $\Xcal\cap B_\varepsilon$, hence compact. 

We show now that the map $s_|$ is Thom regular, see \cref{def:Thom}. We use on $F_f$ the Whitney stratification $\{S_i^\circ,S_i^\partial\}$ introduced before Lemma \ref{lem:lagrange}.

Moreover, $s_|$ is a stratified submersion if we take the stratification $\big\{D_\eta^*,\{0\}\big\}$ on $D_\eta$. Indeed, the stratification of $\Xcal$ is adapted to $X$, so $s_|^{-1}(0)$ is union of strata (including on the boundary). The stratified critical points of $s_|$ in the interior coincide (away from $X$) with intersection points $\Gamma\cap F_f$, but we took $\eta$ small to avoid $\Delta$ and, hence, $\Gamma$. For the boundary, the fibres of $s$ are transverse to the Milnor sphere $S_\varepsilon$ (Remark \ref{rem:C1-consequences}), so the restriction of $s_|$ to strata on the sphere is also submersive. 

For the Thom condition, we consider a pair of strata of $F_f\cap s^{-1}(D_\eta)$, the smaller one being $S_j^\circ$ or $S_j^\partial$, a point $p$ in it and a sequence $\{q_n\}_n$ in the larger one converging to $p$. If the smaller stratum is $S_j^\circ$, then the larger one is also interior, and the pair is a pair of strata of $\Xcal\cap f^{-1}(t_0)$ with its Whitney stratification, which is adapted to $X\cap f^{-1}(t_0)$; hence the condition $a_s$ holds by \cite[Theorem 4.2.1]{Briancon1994} (a local result, which we apply at the point $p$). If the smaller stratum is $S_j^\partial$ and $S_j\not\subseteq X$, then $s|_{S_j^\partial}$ is a submersion at $p$ (there are no tangency points over $D_\eta$) and the Thom condition follows from Whitney's condition (a): $\lim_n\ker D_{q_n}s=\ker(d_ps|_{T})$, where $T$ is the limit of the tangent spaces of the larger stratum, and $T\supseteq T_pS_j^\partial$. Finally, if the smaller stratum is $S_j^\partial$ with $S_j\subseteq X$, so that $\ker D_ps|_{S_j^\partial}=T_pS_j^\partial$, let $\hat S_i:=S_i\cap f^{-1}(t_0)$ be the stratum of $\Xcal\cap f^{-1}(t_0)$ containing the $q_n$, and $\hat S_j:=S_j\cap f^{-1}(t_0)$. By \cite[Theorem 4.2.1]{Briancon1994} applied to the germ of $s$ on $\Xcal\cap f^{-1}(t_0)$ at $p$,
\[ \lim_n \ker D_{q_n}\, s|_{\hat S_i}\supseteq T_p\hat S_j.\]
If the larger stratum is $S_i^\circ$, this is the required condition, since $T_p\hat S_j\supseteq T_pS_j^\partial$. If the larger stratum is $S_i^\partial$, then $\ker D_{q_n}s|_{S_i^\partial}=\ker D_{q_n}s|_{\hat S_i}\cap T_{q_n}S_\varepsilon$ and, by the transversality of the Milnor sphere $S_\varepsilon$ to the fibres of $s$ (Remark \ref{rem:C1-consequences}) and to $\hat S_j$,
\[\begin{aligned}
 \ker D_p\, s|_{S_j^\partial} & = T_p\hat S_j \cap T_p S_\varepsilon\\
& \subseteq \left(\lim_n \ker D_{q_n}\, s|_{\hat S_i}\right)\cap T_p S_\varepsilon\\
 & = \lim_n \left(\ker D_{q_n}\, s|_{\hat S_i}\cap T_{q_n} S_\varepsilon\right) = \lim_n \ker D_{q_n}\, s|_{S_i^\partial},
\end{aligned}\]
where the third equality holds because the subspaces $\ker D_{q_n}\, s|_{\hat S_i}$ and their limit are transverse to the sphere, so the intersections have constant dimension.

Therefore, $s_|\colon F_f\cap s^{-1}(D_\eta)\to D_\eta$ is a proper stratified Thom map, which is what the lifting of the radial vector field requires.

\begin{figure}[htbp]
\centering
\begin{tikzpicture}[x=.83cm,y=.83cm,line cap=round,line join=round]
 \draw[->,gray] (-.5,0)--(4.6,0) node[right] {$u=f$};
 \draw[->,gray] (0,-2.6)--(0,2.8) node[above] {$v=s$};
 \draw[thick,draw={rgb,255:red,30;green,140;blue,30}] (0,0) .. controls (.55,1.08) and (2.8,1.92) .. (4,1.84);
 \draw[thick,draw={rgb,255:red,30;green,140;blue,30}] (0,0) .. controls (.7,.06) and (2.1,.64) .. (4,1.36);
\draw[thick,draw={rgb,255:red,30;green,140;blue,30}] plot[smooth] coordinates {
  (-.5,2.2)
	(2,2.4)
	(3.5,2.3)
  (4,2.4)
};
\draw[thick,draw={rgb,255:red,30;green,140;blue,30}] (0,0) .. controls (.7,-0.06) and (2.1,-0.64) .. (4,-1.36);
 \draw[thick,draw={rgb,255:red,30;green,140;blue,30}] (0,0) .. controls (.55,-1.08) and (2.8,-1.92) .. (4,-1.84);
 \node[text={rgb,255:red,30;green,140;blue,30}] at (4.28,1.6) {$\Delta$};
 \draw[dashed,thick,draw={rgb,255:red,140;green,30;blue,30}] (2.1,-2.6)--(2.1,2.6);
 \node[below,text={rgb,255:red,140;green,30;blue,30}] at (2.1,-2.6) {$\{u=t_0\}$};
 \fill (2.1,0) circle (1.8pt);
 \draw[thick,draw={rgb,255:red,220;green,50;blue,50}] plot[smooth] coordinates {
  (1.86,2.05)
	(2.1,2.1)
  (2.34,2.05)
};
\node[left,text={rgb,255:red,220;green,50;blue,50}] at (1.86,2) {$D_{\delta'}$};
\draw[thick,draw={rgb,255:red,220;green,50;blue,50}] plot[smooth] coordinates {
  (1.86,-2.05)
	(2.1,-2.1)
  (2.34,-2.05)
};
 \draw[thick,draw={rgb,255:red,220;green,50;blue,50}] plot[smooth] coordinates {
  (1.86,.4)
	(2.1,.45)
  (2.34,.4)
};
\node[right,text={rgb,255:red,220;green,50;blue,50}] at (2.2,.3) {$D_\eta$};
\draw[thick,draw={rgb,255:red,220;green,50;blue,50}] plot[smooth] coordinates {
  (1.86,-.4)
	(2.1,-.45)
  (2.34,-.4)
};
\end{tikzpicture}
\caption{Sketch of the different elements in the proof.}
\label{fig:sketch}
\end{figure}
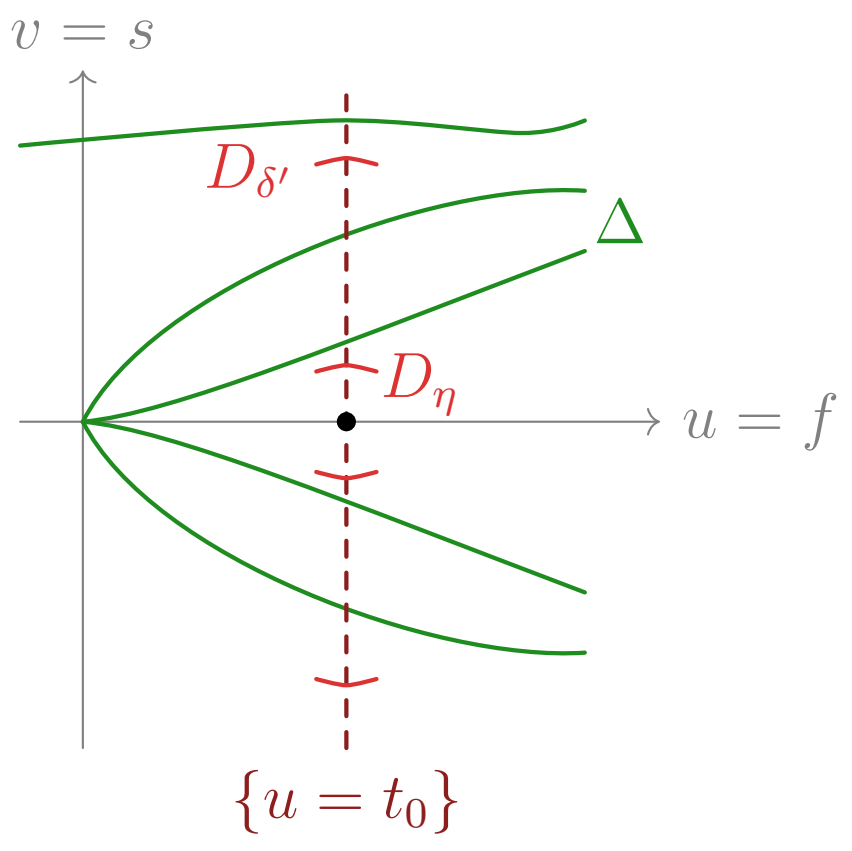

In the last step of the proof we apply stratified Morse theory on $F'_f$. 
We denote by $s_|$ the restriction of $s$ to $F'_f$. Since a holomorphic function cannot be modified only in a neighbourhood of a compact set, we take as Morsification of $s_|$ the function $s^m_|:=(s+\ell)|_{F'_f}$, where $\ell$ is a generic linear form on $\C^N$, as small as needed. For $\ell$ small enough, $s^m_|$ has the following properties:
\begin{enumerate}
\item its critical points on the interior strata are stratified Morse points lying in a small neighbourhood of the finite set $\Gamma\cap F_f$, and $n_i$ of them lie in the stratum $S_i$ (Definition \ref{def:ni-chii});
\item $|s^m_||$ has no tangency points on $\partial F'_f$, and $s^m_|$ is a stratified submersion on the regions $\{\eta-c\le|s|\le\eta+c\}$ and $\{\delta'-c\le|s|\le\delta'\}$, for some small $c>0$; indeed, these are the corresponding properties of $s_|$ (Remark \ref{rem:C1-consequences}), which are open conditions on the compact set $F'_f$.
\end{enumerate}
By (2) and the Thom--Mather first isotopy lemma, $s_|$ is a trivial fibration over the two annuli, so that 
\[
\{|s^m_|| \le \eta\}\simeq\{|s_||\le\eta\}\simeq \{|s_||\le\eta'\} 
\qquad\text{for all } \eta-c\le\eta'\le\eta+c,
\]
since the first set is squeezed between $\{|s_||\le\eta-c\}$ and $\{|s_||\le\eta+c\}$ (for $\ell$ small), and analogously $\{|s^m_||\le\delta'\}\simeq F'_f$.

Now we apply Lemma \ref{lem:no-contribution} to $g^m:=|s^m_||$ on $F'_f$: by (2), condition (2) of the lemma is void, so that between two consecutive critical values of $g^m$ in $[\eta,\delta']$ the sublevel sets have the same homotopy type. At each critical value, the critical points are the stratified Morse points of $s^m_|$ given in (1), and by stratified Morse theory (see, e.g., \cite[Theorem 5.5.3]{Goresky2020}, cf. \cite[Chapter 5]{Goresky1988}) passing through a Morse point in the stratum $S_i$ attaches a stratified handle whose Euler characteristic is $\chi_i$: the tangential Morse data has Euler characteristic $(-1)^{\dim S_i-1}$, because the real index of $|s^m_||$ at a Morse point of the holomorphic function $s^m_|$ on the stratum $S_i\cap f^{-1}(t_0)$, of complex dimension $\dim S_i-1$, is $\dim S_i-1$, and the normal Morse data is the pair $(\operatorname{Cone}(\mathcal L_i),\mathcal L_i)$. By the additivity of the Euler characteristic we get
\begin{equation}\label{eq:handle-comparison}
 \chi\bigg(\big(g^m\big)^{-1}[0,\eta]\bigg)
       +\sum_i\chi_i n_i
 =\chi\bigg(\big(g^m\big)^{-1}[0,\delta']\bigg).
\end{equation}
Using the homotopy equivalences 
\[\begin{aligned}
(g^m)^{-1}[0,\eta] & \simeq   F_f\cap X\quad\text{and}\\
(g^m)^{-1}[0,\delta'] &\simeq F'_f\simeq F_f
\end{aligned}\]
 explained above,
\cref{eq:handle-comparison} yields \cref{eq:main-formula}.
\end{proof}

\begin{remark}
The Thom condition is used only in Lemma \ref{lem:C1}, and only along $\Phi^{-1}(0)$. Hence Theorem \ref{thm:general} holds, with the same proof, for any pair $(f,s)$ such that $\Gamma$ is empty or a curve and the fibres of $\Phi$ are transverse to the Milnor sphere over a small polydisc.
\end{remark}

{ 
\begin{remark}\label{rem:Le-Greuel}
In the recent preprint \cite{LNS2025} a general Lê--Greuel formula is given. 
Their formula \cite[Corollary 5.15]{LNS2025} compares $F_s$ with $F_{f,s}$, the Milnor fibre of $\Phi$, i.e.\ $B_\epsilon\cap f^{-1}(u)\cap s^{-1}(v)$, for a generic small $(u,v)$ in $\C^2$. Moreover, their formula counts critical points of $f$ on $F_s$. In contrast, we compare $F_f$ with $F_f\cap X=B_\epsilon\cap f^{-1}(u)\cap s^{-1}(0)$, with $u\ne0$, and we count critical points of $s$ on $F_f\cap \Gamma$.
\end{remark}

The following example illustrates that our theorem also works even when $\Crit \Phi$ has an irreducible component inside $f^{-1}(0)$.

\begin{example}\label{ex:xy}
Let $f,s\colon(\CC^3,0)\to(\CC,0)$ be such that $f(x,y,z)=xy$ and $s(x,y,z)=z$. We see that $f|_X$ has an isolated critical point, so Theorem \ref{thm:general} applies.  We have $\Gamma=\Crit \Phi=\{x=y=0\}\subset f^{-1}(0)$. 
Moreover, 
\[
F_f\cong (F_f\cap X)\times D,
\] 
but $s$ is regular on the interior of $F_f$, so both sides of \cref{eq:main-formula} vanish.\end{example}

\begin{example}\label{ex:two-components}
Let $f,s\colon(\CC^3,0)\to(\CC,0)$ be given by $f(x,y,z)=x^2+y^2(y+z)$ and $s(x,y,z)=z$. Then $f|_X=x^2+y^3$ has an isolated critical point (an $A_2$ singularity), so Theorem \ref{thm:general} applies. We have 
\[\begin{aligned}
\Crit\Phi&=\{x=0,\ y(3y+2z)=0\}=L_1\cup L_2,\\
 L_1&=\{x=y=0\},\\
 L_2&=\{x=0,\ 3y+2z=0\},
\end{aligned}\]
so $\Gamma=L_1\cup L_2$, where $L_1\subset f^{-1}(0)$, while $f=4z^3/27$ on $L_2$. After the change of coordinates $w=y+z$, $f=x^2+y^2w$ is the $D_\infty$ singularity, whose Milnor fibre is a disc bundle over $S^2$, so $\chi(F_f)=2$; on the other hand, $F_f\cap X$ is the Milnor fibre of $x^2+y^3$, so $\chi(F_f\cap X)=1-2=-1$. Hence the left hand side of \eqref{eq:main-formula} is $3$. The critical points of $s|_{F_f}$ are the three points of $L_2\cap f^{-1}(t)$, given by $4z^3=27t$, each of them a Morse point, so $n_0=3$ for the open stratum $S_0$, $\chi_0=(-1)^{3-1}=1$, and the right hand side is also $3$. We will see in Example \ref{ex:two-components-bis} how to obtain this number algebraically.
\end{example}

\section{Milnor spaces}

Our general theorem applies for a kind of spaces that have a property called the \textit{Milnor condition}. In that setting, the theorem takes a very easy form that we explain below.
\newline

A \emph{Milnor space} is a reduced complex space $\Xcal$ such that at each point $x\in \Xcal$, the rectified homotopical 
depth $\rhd(\Xcal,x)$ is equal to $\dim(\Xcal,x)$. A general reference on this topic is \cite{Hamm-Le-Handbook}, where the definition of the rectified homotopical depth and basic properties of Milnor spaces are explained. 

In particular, examples of Milnor spaces are smooth spaces and hypersurfaces in Milnor spaces (defined locally by one equation). Hence, any local complete intersection, with arbitrary singularities, is a Milnor space.

\begin{theorem}[see {\cite[Theorem 9.5.4]{Hamm-Le-Handbook}}]\label{Hamm-Le} 
Let $(\Xcal,x)$ be a germ of Milnor space of dimension $n+1$ and assume that $f\colon(\Xcal,x)\to(\mathbb C,0)$ has an isolated critical point in the stratified sense. Then, the general fibre $F_f$ of $f$ has the homotopy type of a bouquet of spheres of dimension $n$.
\end{theorem}

\begin{definition} 
With the hypothesis and notation of \cref{Hamm-Le}, the number of spheres $\beta_n(F_f)$ is called the \emph{Milnor number} of $f$ and is denoted by $\mu(f)$.
\end{definition}

Hence, the Euler characteristic $\chi(F_f)$ is equal to $1+(-1)^n\mu(f)$, which simplifies the left hand side of \cref{eq:main-formula}. For the right hand side we also have one simplification because the 
 {complex link $\mathcal L_i$ of the stratum $S_i$ is} homotopically a wedge of spheres of dimension $\dim \Xcal -\dim S_i -1$, see \cite[Theorem 9.5.1]{Hamm-Le-Handbook}.
 {If $d_i:=\dim S_i$ and $\beta_i:=\beta_{n-d_i}(\mathcal L_i)$,
\begin{align*}
\chi_i&=(-1)^{d_i-1}
     \chi\bigl(\operatorname{Cone}(\mathcal L_i),\mathcal L_i\bigr)\\
     &=(-1)^{d_i-1}\left(1-(1+(-1)^{n-d_i}\beta_i) \right)\\
     &=(-1)^n \beta_i.
   \end{align*}   
}
For the open strata (where $d_i=n+1$ and $\mathcal L_i=\varnothing=S^{-1}$) we set $\beta_i:=1$, which is consistent with the above computation, since $\chi(\operatorname{Cone}(\varnothing),\varnothing)=1$.

{ 
\begin{corollary}\label{cor:Milnor}
Let $(\Xcal,0)$ be a Milnor space of dimension $n+1$, let $\mathcal S$ be a Whitney stratification of $\Xcal$ adapted to the hypersurface $X=s^{-1}(0)$, and let $f,s$ be functions on $\Xcal$ such that $f$ and $f|_X$ have an isolated (stratified) critical point with respect to $\mathcal S$ (the latter being automatic, since $\mathcal S$ is adapted to $X$). Then
\[\mu(f)+\mu(f|_X)=\sum_i n_i\beta_i.\]
\end{corollary}

We remark that in Corollary \ref{cor:Milnor}, all the terms that appear in the formula are non-negative. This is important when one wants to study properties like upper semicontinuity of the invariants.
}

\section{An algebraic formula}

Throughout this section, we will assume that $\Xcal$ has pure dimension and that its singular locus $\Sing (\Xcal)$ is contained in $f^{-1}(0)$. Under these assumptions, the fibre $F_f$ is a smooth manifold with boundary. We will provide an algebraic formula to compute the difference $\chi(F_f)-\chi(F_{f}\cap X)$ in terms of the Jacobian ideal of the pair $(f,s)$ on $\Xcal$. Such a formula is quite useful because it can be implemented easily in any computer algebra system like \textsc{Singular} \cite{DGPS}.

In order to do it, we need to recall a couple of definitions. The first one is the Jacobian ideal $J_{\Xcal}(h)$ of a map germ $h\colon(\Xcal,0)\to(\C^r,0)$. 
Suppose that $(\Xcal,0)\subseteq(\C^N,0)$ and that $\phi_1,\dots,\phi_k\in \Oscr_N$ generate the 
ideal of vanishing functions $I_{\Xcal,0}$ in $\Oscr_N$ and let $\phi=(\phi_1,\dots,\phi_k)$. Take also a holomorphic extension $\bar h\colon (\C^N,0)\to(\C^r,0)$ of $h$. Then \textit{the Jacobian ideal} $J_{\Xcal}(h)\subseteq \mathcal O_{\Xcal,0}$ is given by
\[
J_{\Xcal}(h)=\frac{I_{\Xcal,0}+J_{N-d+r}(\phi,\bar h)}{I_{\Xcal,0}},
\]
where $d=\dim(\Xcal,0)$ and $J_{N-d+r}(\phi,\bar h)$ is the ideal in $\mathcal O_N$ generated by the $(N-d+r)$-minors of the Jacobian matrix of $(\phi,\bar h)$. 

It is not difficult to check that the definition of $J_{\Xcal}(h)$ is independent of the choice of the generators $\phi_1,\dots,\phi_k$ and the extension $\bar{h}$. Moreover, 
by the Jacobian criterion, the zero locus $V(J_{\Xcal}(h))$ is the set germ of points $x\in \Xcal$ such that either $x\in \Sing (\Xcal)$ or $x\notin \Sing (\Xcal)$ and $x$ is a critical point of
$h\colon \Xcal\to \C^r$.

The second definition we recall is the notion of saturation of $I$ with respect to $J$, where $I,J$ are ideals in a ring $R$. It is the ideal 
\[
I\colon J^\infty=\left\{a\in R\ |\ aJ^k\subseteq I, \text{ for some $k\ge 1$}\right\}.
\]
When $R=\mathcal O_{\Xcal,0}$ we have a nice geometrical interpretation: 
$$V(I\colon J^\infty)=\overline{V(I)\setminus V(J)}\,,$$
where the bar means the Zariski closure (see e.g. \cite{Greuel2008}).

\begin{theorem}\label{thm:algebraic} Suppose that $(\Xcal,0)$ has pure dimension $n+1$ and let $f,s\colon (\Xcal,0)\to(\C,0)$ be two holomorphic functions such that 
$\Sing(\Xcal)\subseteq f^{-1}(0)$ and $f|_X$ has an isolated critical point, where $X=s^{-1}(0)$.
Then,
\[
\chi(F_f)-\chi(F_{f}\cap X)=(-1)^n \dim_\C\frac{\mathcal O_{\Xcal,0}}{(f)+J_{\Xcal}(f,s)\colon (f\cdot s)^\infty} \;.
\]
\end{theorem}

\begin{proof}
Since $\Xcal$ has pure dimension $n+1$ and $\Sing (\Xcal)\subseteq f^{-1}(0)$, we can assume that $\Xcal\setminus f^{-1}(0)\cup X$ is a single stratum $S_0$ of the refined stratification $\mathcal S'$ of $\Xcal$ (equivalently, $\Xcal\setminus X$ is a single stratum of $\mathcal S$). 
Since $F_f\subseteq \Xcal\setminus f^{-1}(0)$ and the Morse points of a Morsification of $s$ on $F_f$ lie near $\Gamma\cap F_f$, away from $X$,
Theorem \ref{thm:general} gives
\[
\chi(F_f)-\chi(F_{f}\cap X)=(-1)^n n_0,
\]
where $n_0$ is the number of critical points of a Morsification of $s$ on $F_f$ near the finite set $F_f\cap \Gamma$. By the conservation of the Milnor number, $n_0$ 
is equal to the sum of the Milnor numbers of the critical points of $s$ on $F_f\cap \Gamma$. More precisely, suppose that $F_f=f^{-1}(u)\cap B_\epsilon$ for some $u\ne0$ and put $Z_u=f^{-1}(u)$. Then, 
\begin{equation}\label{eq:sum}
n_0=\sum_{x\in F_f\cap \Gamma} \mu(s|_{Z_u};x)=\sum_{x\in F_f\cap \Gamma}\dim_\C\frac{\mathcal O_{Z_u,x}}{J_{Z_u}(s|_{Z_u})}.
\end{equation}
Moreover, $\Xcal$ is smooth at any point $x\in Z_u$, which implies that
\begin{equation}\label{eq:milnor}
\frac{\mathcal O_{Z_u,x}}{J_{Z_u}(s|_{Z_u})}\cong \frac{\mathcal O_{\Xcal ,x}}{(f-u)+J_{\Xcal} (f,s)}.
\end{equation}

By construction, the zero locus of the ideal $J_{\Xcal}(f,s)\colon (f\cdot s)^{\infty}$ is 
\[
V(J_{\Xcal}(f,s)\colon (f\cdot s)^{\infty})=\overline{V(J_{\Xcal}(f,s))\setminus V(f\cdot s)}=\Gamma',
\]
where $\Gamma':=\overline{\Gamma\setminus f^{-1}(0)}$ is the union of the irreducible components of $\Gamma$ not contained in $f^{-1}(0)$; here we use that $\Sing(\Xcal)\subseteq V(f)$ and that $V(J_{\Xcal}(f,s))\setminus\Sing(\Xcal)=\Crit\Phi\cap S_0$. Note that $\Gamma'\cap F_f=\Gamma\cap F_f$. We consider $\Gamma'$ as a complex space with the analytic structure induced by $J_{\Xcal}(f,s)\colon (f\cdot s)^{\infty}$, so that its local ring is $\mathcal O_{\Gamma'}=\mathcal O_{\Xcal,0}/J_{\Xcal}(f,s)\colon (f\cdot s)^{\infty}$.

By Theorem \ref{thm:LNS}, $\Gamma'$ is empty or a curve. If $\Gamma'=\varnothing$, both sides of the formula vanish, so we assume that $\Gamma'$ is a curve. Hence all the primary components of $J_{\Xcal}(f,s)\colon (f\cdot s)^{\infty}$ in $\mathcal O_{\Xcal,0}$ have dimension $\le 1$. Assume that $J_{\Xcal}(f,s)$ has an irredundant primary decomposition in $\mathcal O_{\Xcal,0}$ as follows:
\[
J_{\Xcal}(f,s)=\mathfrak q_1\cap\dots\cap\mathfrak q_t.
\]
It is well known that the saturation commutes with finite intersections, which gives
\[
J_{\Xcal}(f,s)\colon (f\cdot s)^{\infty}=\bigcap_{i} (\mathfrak q_i\colon (f\cdot s)^{\infty}).
\]
A simple computation shows that $\mathfrak q_i\colon (f\cdot s)^{\infty}=\mathcal O_{\Xcal,0}$ when $f\cdot s\in \sqrt{\mathfrak q_i}$ and $\mathfrak q_i\colon (f\cdot s)^{\infty}=\mathfrak q_i$ otherwise.
In particular, the saturation kills all the $\mathfrak m$-primary 
components and thus, $J_{\Xcal}(f,s)\colon (f\cdot s)^{\infty}$ does not have any $\mathfrak m$-primary 
component, that is, all its primary components have dimension $1$. This shows that  $\mathcal O_{\Gamma'}$ is unmixed of dimension 1. 

Moreover, since $\Gamma'\cap V(f)=\{0\}$ (each irreducible component of $\Gamma'$ is a curve not contained in $V(f)$), $f$ induces a regular element $\bar{f}$ in $\mathcal O_{\Gamma'}$, and thus $\mathcal O_{\Gamma'}$ is Cohen--Macaulay (see \cite[Exercise C.3.1]{Mond2020}). Finally, we use the conservation of the multiplicity (see \cite[Corollary E.6]{Mond2020}):
 \begin{align*}
\dim_\C\frac{\mathcal O_{\Xcal,0}}{(f)+J_{\Xcal}(f,s)\colon (f\cdot s)^\infty}&=
\dim_\C\frac{\mathcal O_{\Gamma'}}{(\bar f)}\\
&=\sum_{x\in \bar f^{-1}(u)\cap{B}_\epsilon}\dim_\C\frac{\mathcal O_{\Gamma',x}}{(\bar f-u)}\\
&=\sum_{x\in F_f\cap\Gamma}\dim_\C\frac{\mathcal O_{\Xcal,x}}{(f-u)+J_{\Xcal}(f,s)\colon (f\cdot s)^\infty}\\
&=\sum_{x\in F_f\cap\Gamma}\dim_\C\frac{\mathcal O_{\Xcal,x}}{(f-u)+J_{\Xcal}(f,s)}=n_0,
\end{align*}
by equations \eqref{eq:sum} and \eqref{eq:milnor}.
The fourth equality holds because $f\cdot s$ is a unit in the ring $\mathcal O_{\Xcal,x}$ when $x\in F_f\cap\Gamma$. 

\end{proof}

\begin{remark} 
It follows from the proof of \cref{thm:algebraic} that the dimension on the right hand side of the formula has a nice geometrical interpretation as a local intersection number:
\[
i(V(f),\Gamma';0)=\dim_\C\frac{\mathcal O_{\Xcal,0}}{(f)+J_{\Xcal}(f,s)\colon (f\cdot s)^\infty}.
\]
Here  $i(V(f),\Gamma';0)$ is the local intersection number at the origin in $\Xcal$ of the hypersurface $V(f)$ and the  curve $\Gamma'$, when considered as a complex space with the analytic structure induced by $J_{\Xcal}(f,s)\colon (f\cdot s)^{\infty}$. Now the additivity of the local intersection number applies: if $\Gamma'$ has irreducible components parametrised by curves $\gamma_i\colon(\C,0)\to(\Xcal,0)$, for $i=1,\dots,r$, then
\[
i(V(f),\Gamma';0)=\sum_{i=1}^r m_i \ord(f\circ \gamma_i),
\]
where $m_i$ is the multiplicity of $\Gamma'$ along $\gamma_i$ and $\ord(f\circ \gamma_i)$ is the order of $f\circ \gamma_i$.
\end{remark}

\begin{example}\label{ex:two-components-bis}
In Example \ref{ex:two-components}, 
\[\begin{aligned} 
J_{\Xcal}(f,s)=&\big(x,\,y(3y+2z)\big),\quad\text{and}\\
J_{\Xcal}(f,s)\colon (f\cdot s)^\infty=&(x,\,3y+2z).\end{aligned}\]
The saturation removes the component $L_1\subset f^{-1}(0)$, so that $\Gamma'=L_2$ with its reduced structure. Theorem \ref{thm:algebraic} gives
\[
\chi(F_f)-\chi(F_f\cap X)=\dim_\C\frac{\C\{y,z\}}{(y^2(y+z),\,3y+2z)}=\dim_\C\frac{\C\{z\}}{(z^3)}=3,
\]
in agreement with Example \ref{ex:two-components}. Note that the saturation with respect to $s$ alone does not remove $L_1$. Indeed, $J_{\Xcal}(f,s)\colon s^\infty=J_{\Xcal}(f,s)$, and $\mathcal O_{\Xcal,0}/\big((f)+J_{\Xcal}(f,s)\big)$ is infinite dimensional. This is why the saturation with respect to $f\cdot s$ is needed in Theorem \ref{thm:algebraic} when $f$ does not have an isolated critical point (cf.\ Lemma \ref{lem:Gamma} below).
\end{example}

As in Section 3, we can apply \cref{thm:algebraic} in the case that $(\Xcal,0)$ is a Milnor space and  both $f$ and $f|_X$ have an isolated critical point in the stratified sense (see \cref{cor:Milnor}). We remark that a Milnor space $(\Xcal,0)$ is pure dimensional, so it is not necessary to add such condition (see \cite[Subsection 9.5.1]{Hamm-Le-Handbook}). 

In fact, when $f$ has an isolated critical point, the hypothesis $\Sing (\Xcal)\subseteq f^{-1}(0)$ is equivalent to $\Xcal$ having an isolated singularity, since all singular points of $\Xcal$ are critical points of $f$ (the strata of $\mathcal S$ contained in $\Sing(\Xcal)\cap f^{-1}(0)$ would be contained in $\Crit f$).

Moreover, we can use saturation with respect to $s$ instead of $f\cdot s$, as the following lemma shows.

\begin{lemma}\label{lem:Gamma} Suppose that $\Sing(\Xcal)\subseteq f^{-1}(0)$ and that $f$ has an isolated critical point (with respect to $\mathcal S$). Then, we have the equality $\Gamma=\Gamma'$, where
 \[
 \Gamma'=\overline{\Gamma\setminus f^{-1}(0)}.
\]
Moreover, we also have $J_{\Xcal}(f,s)\colon (f\cdot s)^\infty=J_{\Xcal}(f,s)\colon s^\infty$.
\end{lemma}

\begin{proof}
The inclusion $\Gamma'\subseteq \Gamma$ is obvious. Suppose, by contradiction, that $\Gamma'\subsetneq\Gamma$. By the curve selection lemma, we could find an analytic curve $\gamma:(\C,0)\to(\Xcal,0)$ such that $\gamma(t)\in \Crit\Phi\setminus X$ and $\gamma(t)\notin\Gamma'$, for all $t\ne0$. 
Since $\Crit\Phi\setminus (X\cup f^{-1}(0))\subseteq\Gamma\setminus f^{-1}(0)\subseteq\Gamma'$, we must have $f(\gamma(t))=0$ and hence, $D_{\gamma(t)} f(\gamma'(t))=0$, 
by the chain rule. 

We can assume, without loss of generality, that there exists a single stratum $S_i$ such that $\gamma(t)\in S_i$, for all $t\ne0$. By definition, $\gamma(t)$ is a critical point of $\Phi|_{S_i}$, so $D_{\gamma(t)} (f|_{S_i})$ and $D_{\gamma(t)} (s|_{S_i})$ are collinear. Since $f$ has an isolated critical point, $D_{\gamma(t)} (f|_{S_i})\ne0$ and therefore, 
\[
D_{\gamma(t)} (s|_{S_i})=\lambda(t) D_{\gamma(t)} (f|_{S_i}),
\]
for some $\lambda(t)\in\C$. We deduce that $D_{\gamma(t)} s(\gamma'(t))=0$, hence $s(\gamma(t))=0$ and $\gamma(t)\in X$, in contradiction with the choice of $\gamma$.

For the equality of the ideals, let $J_{\Xcal}(f,s)=\mathfrak q_1\cap\dots\cap\mathfrak q_t$ be an irredundant primary decomposition. As in the proof of Theorem \ref{thm:algebraic}, $\mathfrak q_i$ survives in $J_{\Xcal}(f,s)\colon (f\cdot s)^\infty$ if and only if $f\cdot s\notin\sqrt{\mathfrak q_i}$, and in $J_{\Xcal}(f,s)\colon s^\infty$ if and only if $s\notin\sqrt{\mathfrak q_i}$. The two saturations can only differ by the components $\mathfrak q_i$ with $V(\mathfrak q_i)\subseteq V(f)$ and $V(\mathfrak q_i)\not\subseteq X$. Such a $V(\mathfrak q_i)$ is contained in an irreducible component $C$ of $V(J_{\Xcal}(f,s))=\Sing(\Xcal)\cup\Crit(\Phi|_{\Xcal\setminus\Sing(\Xcal)})$ with $C\not\subseteq X$. By the previous paragraph, $\Sing(\Xcal)=\{0\}$, so $C\subseteq\overline{\Crit\Phi\setminus X}=\Gamma=\Gamma'$, and then $V(\mathfrak q_i)\subseteq C\cap V(f)\subseteq\Gamma'\cap V(f)=\{0\}$. Hence $\mathfrak q_i$ is $\mathfrak m$-primary and it is killed by both saturations. 
\end{proof}
}

Thus, we have:

\begin{corollary}\label{cor:algebraic} Suppose that $(\Xcal,0)$ is a Milnor space of dimension $n+1$ with an isolated singularity and let $f,s\colon (\Xcal,0)\to(\C,0)$ be two holomorphic functions such that 
both $f$ and $f|_X$ have an isolated critical point, where $X=s^{-1}(0)$. Then,
\[
\mu(f)+\mu(f|_X)=\dim_\C\frac{\mathcal O_{\Xcal,0}}{(f)+J_{\Xcal}(f,s)\colon s^\infty} \;.
\]
\end{corollary}

\section{Some applications}

In this last section, we apply our results to recover Siersma's formula and also to show how to compute the image Milnor number $\mu_I(f)$ of a map germ with isolated instability.

\subsection{Siersma's formula}

In this subsection we consider $\Xcal=\C^{n+1}\times\C$ and the function $f$ is the projection $\pi\colon (\C^{n+1}\times\C,0)\to(\C,0)$ given by $\pi(x,t)=t$. Now take a holomorphic function $s\colon(\C^{n+1}\times\C,0)\to(\C,0)$ and consider the hypersurface $X=s^{-1}(0)$. 

We choose a Whitney stratification of $X$ and assume that $\pi|_X$ has an isolated critical point in the stratified sense. This implies that the Milnor fibre $F_{\pi|_{X}}$ has the homotopy type of a wedge of $n$-spheres and the Milnor number $\mu(\pi |_X)$ is the number of such spheres, $\beta_n(B_\epsilon\cap X_t)$, as we saw in Section 3.

However, in this case, we can see $s$ as a family of functions $\{s_t\}$ on $\C^{n+1}$, where $s_t(x)=s(x,t)$ and we can also see $X$ as a family of hypersurfaces $\{X_t\}$ in $\C^{n+1}$, where $X_t=s_t^{-1}(0)$. With this notation, $F_{\pi|_{X}}$ is equal to $B_\epsilon\cap X_t$, for $0<|t|<\delta$ and $0<\delta\ll\epsilon\ll 1$.

\begin{corollary}\label{cor:Siersma} Let $s\colon(\C^{n+1}\times\C,0)\to(\C,0)$ be a holomorphic function such that $\pi |_X$ has an isolated critical point, where $X=s^{-1}(0)$. Then, $\beta_n(B_\epsilon\cap X_t)$ is equal to the sum of the Milnor numbers of the critical points of $s_t $ on $B_\epsilon\setminus X_t$. Moreover, 
\[
\beta_n(B_\epsilon\cap X_t)=\dim_\C\frac{\mathcal O_{n+2}}{(t)+J_x(s)\colon s^\infty} \;,
\]
where $J_x(s)$ is the ideal in $\mathcal O_{n+2}$ generated by the partial derivatives of $s$ with respect to $x_1,\dots,x_{n+1}$. 
\end{corollary}

\subsection{The image Milnor number}

In this subsection we consider singularities of holomorphic map germs $f\colon(\C^n,S)\to(\C^{p},0)$, where $S\subset \C^n$ is a finite subset, in the context of the Thom--Mather theory. We refer to the book \cite{Mond2020} for precise definitions and basic properties on $\Ascr$-equivalence, stability, versality, etc. 

Roughly speaking two map germs are $\Ascr$-equivalent if they are equal up to holomorphic changes of coordinates in the source and the target. A map germ is stable if any unfolding is $\Ascr$-trivial, so stable map germs play the role of rigid objects in the deformation theory of singularities of mappings. When $f$ is unstable, the idea is to perturb it in order to get another mapping which is close enough to $f$ and has only stable singularities. 
This stable perturbation exists and is well defined when $f$ has isolated instability and $(n,p)$ are nice dimensions in the sense of Mather \cite{Mather1971}. The geometrical invariants which appear in the stable perturbation are in fact powerful invariants of the map germ $f$ itself.

The image Milnor number was introduced by Mond \cite{Mond1991} for a 
map germ $f\colon(\C^n,S)\to(\C^{n+1},0)$ with isolated instability with $n<15$ (i.e., 
$(n,n+1)$ are nice dimensions). A stabilisation is a 1-parameter unfolding $F\colon(\C^n\times\C,S)\to(\C^{n+1}\times\C,0)$ given by $F(x,t)=(f_t(x),t)$ with the property that $f_0=f$ and $f_t$ has only stable singularities when $t\ne0$. 

Based on the works of L\^e and Siersma, Mond proved in \cite{Mond1991} that $X_t\cap B_\epsilon $ has the homotopy type of a bouquet of $n$-spheres, where $X_t$ is the image of $f_t$, for $0<|t|<\delta$ and $0<\delta\ll\epsilon\ll 1$. The number of spheres $\beta_n(X_t\cap B_\epsilon)$ is called the image Milnor number and is denoted by $\mu_I(f)$. Mond showed that $\mu_I(f)$ is well defined, that is, it is independent of the choice of $t,\delta,\epsilon$ and the stabilisation $F$. In some sense, the image Milnor number can be considered as a Milnor number of the hypersurface parametrised by the image of $f$, but when we deform the parametrisation instead of deforming the implicit equation.

Moreover, there is a long standing conjecture by Mond \cite{Mond1991} which claims that 
\[
\codim_{\mathscr A_e}(f)\le \mu_I(f),
\]
with equality when $f$ is weighted homogeneous. This is an inequality of type $\tau\le \mu$ for singularities of mappings. Here, $\codim_{\mathscr A_e}(f)$  is the extended $\mathscr A$-codimension and can be interpreted geometrically as the number of parameters in a miniversal unfolding of $f$. Hence, $\codim_{\mathscr A_e}(f)$ is analogous to the Tjurina number $\tau$ of a hypersurface with isolated singularity. Mond's conjecture is known to be true for $n=1,2$ but it is still open for $n\ge 3$. We refer to the survey \cite{Nuno2024} for a recent account on Mond's conjecture.

We use Siersma's formula to compute algebraically the image Milnor number as follows:
denote by $(X,0)$ the image of the stabilisation $F$ and by 
$\pi\colon(\C^{n+1}\times\C,0)\to(\C,0)$ the projection onto the parameter $\pi(y,t)=t$. 
With this notation, $X_t\cap B_\epsilon$ is nothing but  the Milnor fibre $F_{\pi|_X}$ of the restriction $\pi|_X$. We also have a natural 
stratification of $X$ given by the stable types on $X\setminus\{0\}$, together with the stratum $\{0\}$ itself (provided that $f$ is unstable). 
The facts that $f$ has isolated instability and $F$ is a stabilisation imply that the stratification is analytically trivial on  $X\setminus\{0\}$ and the strata are transverse to the hyperplane $t=0$. Therefore, $\pi|_X$ has an isolated critical point in the stratified sense, so we can use \cref{cor:Siersma}.

\begin{corollary}\label{cor:mu_I(f)} 
Let $f\colon(\C^n,S)\to(\C^{n+1},0)$ with isolated instability and let $s\colon(\C^{n+1}\times\C,0)\to(\C,0)$ be a reduced equation of the image of a stabilisation $F$ of $f$. Then, $\mu_I(f)$ is equal to the sum of the Milnor numbers of the critical points of $s_t $ on $B_\epsilon\setminus s_t^{-1}(0)$, where $s_t(y)=s(y,t)$. Moreover, 
\[
\mu_I(f)=\dim_\C\frac{\mathcal O_{n+2}}{(t)+J_y(s)\colon s^\infty} \;,
\]
where $J_y(s)$ is the ideal in $\mathcal O_{n+2}$ generated by the partial derivatives of $s$ with respect to $y_1,\dots,y_{n+1}$. 
\end{corollary}

\begin{example} Consider the map germ $f\colon(\C^2,0)\to(\C^3,0)$ given by
$f(x,y)=(x,y^3,xy+y^5)$. This is Mond's $H_2$-singularity and has $\Ascr_e$-codimension 2 (see \cite{Mond1985}). According to \cite{Marar1996}, a stabilisation $F$ of $f$ is obtained as
\[
F(x,y,t)=(x,y^3-ty,xy+y^5-ty^3).
\]
The stable perturbation has two Whitney umbrellas and one triple point. Moreover, its image has the homotopy type of $S^2\vee S^2$, so that $\mu_I(f)=2$, as expected since $f$ is weighted homogeneous (see \cref{fig:H2}):

\begin{figure}[ht]
\begin{center}
\includegraphics[width=6cm]{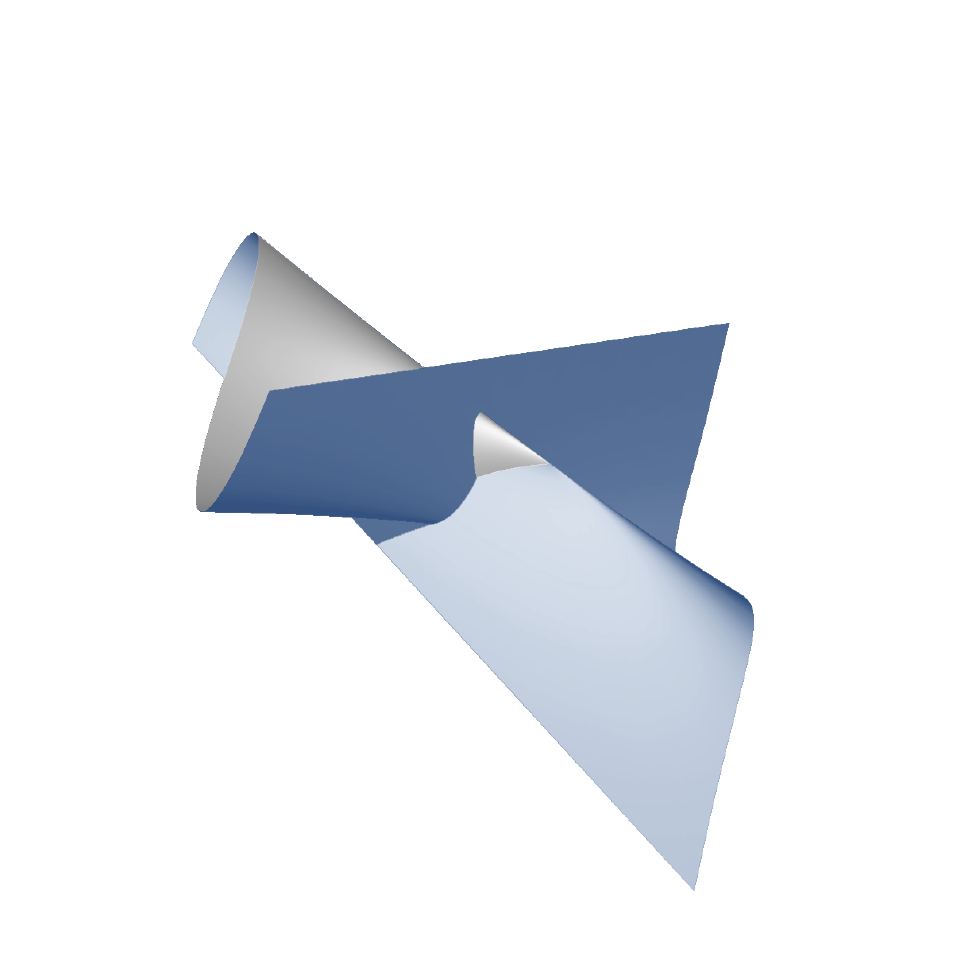}
\caption{Stable perturbation of the $H_2$-singularity}\label{fig:H2}
\end{center}
\end{figure}

We check the formula of Corollary \ref{cor:mu_I(f)}  with the aid of \textsc{Singular}. We denote by $X,Y,Z,T$ the coordinates in the target. Then, the reduced equation $s(X,Y,Z,T)$ of the image of $F$ is
\[
Z^3-X^3Y-3XY^2Z-X^2ZT-2YZ^2T-Y^5+XY^3T+Y^2ZT^2.
\]
Now we compute the relative polar curve $\Gamma$, defined by the ideal $J_y(s)\colon s^\infty$, which has three generators:
\[
X, 3Z-YT, 45Y^2-4T^3.
\]
In fact, in this example $\Gamma$ is an ICIS and the local intersection number of $\Gamma$ with the hyperplane $T=0$ gives $2$, which coincides with $\mu_I(f)$.
\end{example}

\section*{Use of AI tools}
During the preparation of this work, the authors used AI-based language models to assist with checking and drafting parts of the arguments, editing the text, and formatting. The authors take full responsibility for the correctness and content of this paper. AI tools are not credited as authors.

\IfFileExists{myalpha.bst}{\bibliographystyle{myalpha}}{\bibliographystyle{alpha}}
\bibliography{FullBib}
\end{document}

%% file: Style-and-more.tex
\def\dim{\operatorname{dim}}

\def\codim{\operatorname{codim}}

\def\ord{\operatorname{ord}}

\def\rhd{\operatorname{rhd}}

\newcommand{\CC}{\mathbb{C}}

\newcommand{\R}{\mathbb{R}}

\newcommand{\C}{\mathbb{C}}

\newcommand{\Ascr}{\mathscr{A}}

\newcommand{\Oscr}{\mathscr{O}}

\newcommand{\Scal}{\mathcal{S}}

\newcommand{\Xcal}{\mathcal{X}}

\makeatletter
\newcommand{\tpitchfork}{%
  \vbox{
    \baselineskip\z@skip
    \lineskip-.52ex
    \lineskiplimit\maxdimen
    \m@th
    \ialign{##\crcr\hidewidth\smash{$-$}\hidewidth\crcr$\pitchfork$\crcr}
  }%
}
\makeatother

\theoremstyle{plain}
\newtheorem{theorem}{Theorem}[section]

\newaliascnt{lemma}{theorem}
\newtheorem{lemma}[lemma]{Lemma}
\aliascntresetthe{lemma}

\newaliascnt{corollary}{theorem}
\newtheorem{corollary}[corollary]{Corollary}
\aliascntresetthe{corollary}

\newaliascnt{proposition}{theorem}

\aliascntresetthe{proposition}

\newaliascnt{question}{theorem}

\aliascntresetthe{question}

\newaliascnt{thm}{theorem}

\aliascntresetthe{thm}

\newaliascnt{lem}{theorem}

\aliascntresetthe{lem}

\newaliascnt{cor}{theorem}

\aliascntresetthe{cor}

\newaliascnt{prop}{theorem}

\aliascntresetthe{prop}
\Crefname{lemma}{Lemma}{Lemmas}
\crefname{lemma}{Lemma}{Lemmas}
\Crefname{proposition}{Proposition}{Propositions}
\crefname{proposition}{Proposition}{Propositions}
\Crefname{corollary}{Corollary}{Corollaries}
\crefname{corollary}{Corollary}{Corollaries}

\newaliascnt{diagram}{equation}
\aliascntresetthe{diagram}
\crefname{diagram}{Diagram}{Diagrams}
\Crefname{diagram}{Diagram}{Diagrams}

\Crefname{prop}{Proposition}{Propositions}
\crefname{prop}{Proposition}{Propositions}
\Crefname{lem}{Lemma}{Lemmas}
\crefname{lem}{Lemma}{Lemmas}
\Crefname{cor}{Corollary}{Corollaries}
\crefname{cor}{Corollary}{Corollaries}

\Crefname{definition}{Definition}{Definitions}
\crefname{definition}{Definition}{Definitions}

\newtheorem*{theorem**}{Theorem\theoremnum}
\newenvironment{theorem*}[1][]{%
  \edef\theoremnum{\if\relax\detokenize{#1}\relax\else~#1\fi}
  \begin{theorem**}
}{%
  \end{theorem**}
} 
\newtheorem*{conjecture**}{Conjecture\theoremnum}
\newenvironment{conjecture*}[1][]{%
  \edef\theoremnum{\if\relax\detokenize{#1}\relax\else~#1\fi}
  \begin{conjecture**}
}{%
  \end{conjecture**}
} 
\newtheorem*{corollary**}{Corollary\theoremnum}
\newenvironment{corollary*}[1][]{%
  \edef\theoremnum{\if\relax\detokenize{#1}\relax\else~#1\fi}
  \begin{corollary**}
}{%
  \end{corollary**}
}

\theoremstyle{definition}
\newtheorem{definition}[theorem]{Definition}
\newtheorem{proposition/definition}[theorem]{Proposition/Definition}

\newtheorem{example}[theorem]{Example}

\newtheorem{remark}[theorem]{Remark}

\Crefname{notation}{Notation}{Notations}

\theoremstyle{remark}